\documentclass[12pt,twoside,a4paper]{article}

\usepackage{amsmath}
\usepackage{a4wide}
\usepackage{graphicx}
\usepackage{dsfont}
\usepackage{amsthm}
\usepackage{amssymb}

\newcommand{\rem}[1]{}
\newcommand{\todo}[1]{}
\newcommand{\HW}[1]{}

\newcommand{\ui}{\mathrm{i}}

\newcommand{\vz}{\mathbf{z}}

\newcommand{\dof}{{f}}

\title{Cotangent bundle reduction and Poincar{\'e}-Birkhoff normal forms}
\author{\"{U}nver \c{C}ift\c{c}i\thanks{uciftci@nku.edu.tr} $^{1,2}$, Holger Waalkens\thanks{h.waalkens@rug.nl} $^2$ and
Henk Broer\thanks{h.w.broer@rug.nl} $^2$ \\[2ex]
$^1$Department of Mathematics\\ Nam\i k Kemal University\\ 59030 Tekirda\u{g}, Turkey\\[1.5ex]
$^2$Johann Bernoulli Institute for Mathematics and Computer Science\\ University of Groningen\\ PO Box 407\\ 
9700 AK Groningen, The Netherlands}

\begin{document}
\maketitle

\abstract{
In this paper we study a systematic and natural construction of canonical coordinates for the reduced space of a cotangent bundle
with a free Lie group action. 
The canonical coordinates enable us to compute Poincar{\'e}-Birkhoff normal forms of relative equilibria using standard algorithms.
The case of simple mechanical systems with symmetries is studied in detail. 
As examples we compute Poincar{\'e}-Birkhoff normal forms for a Lagrangian equilateral triangle configuration of a three-body system
with a Morse-type potential and the stretched-out configuration of a double spherical pendulum.
}

%\date{\today}

\vspace*{1cm}

%\noindent
%PACS numbers:
%45.50 Jf, % Few and many-body systems
%45.20.Jj	%Lagrangian and Hamiltonian mechanics
%02.40.-k \\
\noindent
AMS classification numbers: 70F07, 70G65, 53C80
\noindent

%%%%%%%%%%%%%%%%%%%%%%
\section{Introduction}

The theory of the reduction of Hamiltonian systems with symmetry is well developed \cite{AbrahamMarsden78, Marsden92}. Although it is a
classical subject and goes back to the pioneers of mechanics, a modern theory was only established in the 1970's.
The main idea can be summarized as follows.

Let $P$ be a symplectic manifold with a symmetry group $G$, and let $J:P \rightarrow \mathfrak{g}^*$ be an equivariant \emph{momentum
mapping} with respect to the coadjoint action of $G$ on $\mathfrak{g}^*$, where $\mathfrak{g}^*$ is the dual space of the Lie algebra $\mathfrak{g}$ of $G$.
Then, under some regularity conditions, the reduced space given by the quotient space $P_\mu:=J^{-1}(\mu)/G_\mu$, where $G_\mu$ is the isotropy group of $\mu \in \mathfrak{g}^*$,
is a symplectic manifold.
A $G$-invariant Hamiltonian function on $P$ can  be reduced to a function on $P_\mu$ which generates the reduced dynamics.

The algebraic definition  as a quotient makes it often difficult to explicitly construct the reduced space $P_\mu$ and develop a good intuition for it.
For example, $P_\mu$  is not necessarily a linear space even if $P$ is linear. This is the case, for example, in the $n$-body problem: 
Although the translation reduced space is Euclidean, the reduced space of rotations is  in general not linear  \cite{LittlejohnReinsch97}. 
But as the reduced space is a symplectic manifold it follows from the  Darboux theorem that  one can locally construct canonical coordinates 
so that  the reduced space locally becomes a linear symplectic space. 
Such canonical coordinates are very useful. For example, they form the starting point of standard algorithms for  the computation of the Poincar{\'e}-Birkhoff normal form at an equilibrium point of a Hamiltonian system \cite{Deprit69,AKN88,DragtFinn76,Murdock03,Broer09}. 
A Poincar{\'e}-Birkhoff normal form is a main tool for the construction of center manifolds and the study bifurcations. As an example of the former application we mention the construction of the phase space structures which govern reaction  dynamics induced by saddle type equilibrium points \cite{Uzeretal02,WaalkensBurbanksWigginsb04}. 

The main objective of this paper is  the systematic construction of
canonical coordinates for the reduced space $P_\mu$ in the case where $P$ is a cotangent bundle and the action of $G$ on $P$ is free, and illustrate how these coordinates can be used to compute Poincar{\'e}-Birkhoff normal forms at the relative equilibria, i.e. the equilibria of the reduced system. This has numerous applications. To give  one example we mention the construction of  the phase space structures which govern the reactions in rotating molecules where the reaction dynamics is induced by saddle type relative equilibria \cite{CiftciWaalkens11}.

We note that the computation of canonical coordinates for a reduced space of a symplectic manifold \cite{RobertsWulffLamb02,RobertsScmahStoica06}
in general or for specific cases such as,
e.g., a cotangent bundle \cite{MarsdenScheurle93,KoonMarsden97} and more concretely for  $n$-body systems \cite{LittlejohnReinsch97,IwaiYamaoka05} have at least implicitly been studied in the literature before.
However, for obtaining the nonlinear terms of a Poincar{\'e}-Birkhoff normal form these works have to be put into context, and a systematic study is missing.
Also the work on the computations of Poincar{\'e}-Birkhoff normal forms of symmetry reduced Hamiltonians is mainly restricted to Abelian Lie group actions.
In this paper we present a systematic approach which covers both the Abelian and the non-Abelian case.

In the following we give a brief review of existing  literature related to this paper.  
In order to obtain canonical coordinates on the reduced space of a cotangent bundle with a Lie group action we follow the method given
in \cite{MarsdenScheurle93} and \cite{KoonMarsden97}  which  take a Lagrangian respectively
Poisson reduction point of view. 
A detailed survey on cotangent bundle reduction and its history can be found in \cite{Marsdenetal07}. 
For the special case of the three-body reduction, our main references are \cite{Iwai87c} and \cite{LittlejohnReinsch97} to which we will come back in Sec.~\ref{sec:Red}.
As for the Poincar{\'e}-Birkhoff normal form, one can find a detailed introduction in \cite{MeyerHall}.
But for completeness, we give a brief explanation of the algorithm in Appendix~\ref{algorithm}.
One of the first applications of the normal form theory to reduced spaces of symplectic spaces with a continuous symmetry
can be found in \cite{ChurchillKummerRod83} where the symmetry group is the circle group.
An application to the restricted three-body problem can be found  in \cite{JorbaVillanueva98} and \cite{Gomezetal04}, for instance.
In \cite{ChencinerFejoz08} normal form computations are done at Lagrange points by using a splitting method.
In another recent work \cite{ChierchiaPinzari11} one can find a detailed study of normal form for planetary systems.
Finally, a normal form at a relative equilibrium of a general dynamical system is given in \cite{LambMelbourne07}. A recent review of
normal form theory in dynamical systems can be found in \cite{Broer09}.   

This paper is organized as follows. 
We start with a general review of the action of Lie groups on tangent and cotangent bundles in Sec.~\ref{sec:1}.
This mainly serves to introduce some basic material and settle the notation.
Sec.~\ref{sec:2} comprises the main result of this paper which is a systematic construction of canonical coordinates for the reduced space of a cotangent bundle with
a free action of  a symmetry group. This includes the derivation of the reduced Hamiltonian in canonical coordinates, a detailed discussion of the case of simple
mechanical systems (Sec.~\ref{sec:simplemechsyst}) and of special cases like Abelian symmetry groups and systems with vanishing angular momenta in Sec.~\ref{sec:examples}, and
the Poincar{\'e}-Birkhoff normal form of relative equilibria in Sec.~\ref{sec:normal}.
Section~\ref{sec:Red} contains our first example which consists of the three-body reduction.
We review in this section how to derive a reduced Hamiltonian in canonical coordinates  in a way which does not depend on the choice of a body-fixed reference frame,
i.e. in the language of Littlejohn and Reinsch  \cite{LittlejohnReinsch97} in a gauge independent way. 
In Sec.~\ref{sec:Normal} we consider a Lagrangian equilateral triangle relative equilibrium, we compute a Poincar{\'e}-Birkhoff normal form at such configurations.
The reconstruction of the full dynamics in the three-body case is addressed in Sec.~\ref{sec:Reconstruction}.
 In Sec.~\ref{sec:Pen} we study our second example which is the double spherical pendulum. After obtaining canonical coordinates for
the reduced system, a normal form computation is done at the relative equilibrium given by the so called stretched out solution. 
Conclusions are  given in Sec.~\ref{sec:conclusions}.    

%%%%%%%%%%%%%%%%%%%%%%
\section{Lie group actions on tangent and cotangent bundles}
\label{sec:1}
In this section we recall the symplectic actions of Lie groups on tangent and
cotangent bundles over a configuration space (mainly to introduce some notation). For the details, we refer to \cite{Arnold78,AbrahamMarsden78, Marsden92, OrtegaRatiu04}. 

Let $G$ be a Lie group and let $M$ be a manifold which is called the \emph{configuration space}. Let the map
\begin{eqnarray}
 G \times M  \rightarrow M \\
(g,s) \rightarrow gs
\end{eqnarray}
be a free action of $G$ on $M$. We denote the left-translation which for a fixed $g\in G$, maps $s\in M$ to
$gs$ by $L_g$.\HW{Henk: this notation is standard.}  The  derived maps of $L_g$ are denoted as follows.
For $s \in M$,
$(L_g)_*:T_sM \rightarrow T_{gs}M$ stands for the derivative map of $L_g$, and $(L_g)^*:T^*_{s}M \rightarrow T^*_{g^{-1}s}M$
stands for the pull-back map of $L_g$. Let $\mathfrak{g}$ denote the Lie algebra of $G$. Then for $\zeta \in \mathfrak{g}$,
the corresponding \emph{infinitesimal generator} or \emph{fundamental vector field} $\zeta_M$ at  $s\in M$ is defined by
\begin{equation}
 \zeta _M(s)= \left.\frac{d}{dt} \right|_{t=0} \left( L_{\exp \left( t\zeta \right)} s \right) \,.
\end{equation}
The $G$ orbit through $s\in M$ is given by $Gs=\{g s | \hspace{2mm} g\in G \} \subset M$. 
The fundamental vector fields $\zeta_M$  are  tangent to the orbits $Gs$ for all $s\in M$.  
Moreover, the tangent space $T_s(Gs)$ is spanned by the fundamental vector fields at $s$.

If $M=G$, i.e. the action is the group operation of $G$, then the fundamental vector fields at $g\in G$ are given by
\begin{equation}
 \zeta _G(g)=\left( R_g \right) _* \zeta \,, \label{Lie} 
\end{equation}
 where  $\zeta \in \mathfrak{g}$ and $R_g$ is the right-translation by $g$ \cite{AbrahamMarsden78}.\HW{Henk: this can be found, e.g., in this reference. The use is standard.}  
 
The coadjoint action of $G$ on the dual space $\mathfrak{g}^*$ of its Lie algebra is defined as
\begin{equation}
\langle (Ad_{g^{-1}})^* \, \mu , \zeta \rangle=\langle \mu , Ad _{g^{-1}} \, \zeta \rangle, \label{eq:def_coadjoint_action}
\end{equation}
for $g\in G$, $\mu \in \mathfrak{g}^*$ and $\zeta \in \mathfrak{g}$. Here $\langle \,,\, \rangle$ stands for the pairing  between a co-vector and vector and
\begin{equation}
Ad_g \, \zeta =\left.\frac{d}{dt} \right|_{t=0} \left( g\, ({\exp \left( t\zeta \right)}) \, g^{-1} \right) \label{eq:def_adjoint_action}
\end{equation}
is the adjoint action. 
The action on $M$ can be lifted to $TM$ and $T^*M$ by the derived maps and both of
the lifted actions are free when the action on $M$ is free. The lifted action on the cotangent bundle is symplectic with respect to
the natural symplectic structure on $T^*M$ \cite{Arnold78}, and has a momentum mapping which is defined as follows:
the \emph{momentum mapping} $J:T^*M \to \mathfrak{g}^*$ is given by
\begin{equation}
 \langle J(s,p), \zeta \rangle=\langle p, \zeta_M(s)  \rangle, \label{momentum}
\end{equation} 
for all $(s,p)\in T^*M$. It is
well-known that $J$ is equivariant with respect to
the action on $T^*M$ and the coadjoint action on $\mathfrak{g}^*$ \cite{Arnold78,AbrahamMarsden78}.

If $M$ is a Riemannian manifold with a Riemannian metric $k$ which is invariant under
the action of $G$, then the lifted action on $TM$ is also symplectic with respect to
the symplectic structure induced by the one on $T^*M$ and the corresponding
momentum mapping $ \mathbf{L}$ is defined by 
\begin{equation}
 \langle \mathbf{L}(s,v), \zeta \rangle= v^T \, k \, \zeta_M (s), \label{tan_mom_map}
\end{equation}  
for all $(s,v) \in TM$ \cite{Marsden92}. 

%%%%%%%%%%%%%%%%%%%%%%
\section{Canonical coordinates in cotangent bundle reduction}
\label{sec:2}
In this section we review the cotangent bundle reduction in the orbit reduction scheme of Marle \cite{Marle76, OrtegaRatiu04}
with a coordinate-based approach. We then obtain canonical coordinates for the reduced space.
A more detailed explanation of the notions used in the following subsection can be found for instance in \cite{MarsdenScheurle93,KoonMarsden97}.

%%%%%%%%%%%%%%%%%%%%%%%%
\subsection{Reduction of  the equations of motion}
\label{orbit}

%We retain the notation of the previous section. 
The \emph{shape space} or \emph{internal space} $Q$ is defined as the quotient $M/G$. As we assume that the action of $G$ on $M$ is free it follows from standard theorems that $M/G$ has a manifold structure and 
 $M \rightarrow Q$
is a fibre bundle \cite{AbrahamMarsden78}.  Using the fibre bundle structure 
one can locally obtain a coordinate system on $M$ by choosing a local coordinate system 
on $Q \times G$.
 Let us assume that a point in $Q \times G$ has
coordinates $(q,g)$.  The coordinates $q$ are called \emph{shape coordinates} or \emph{internal coordinates}.
Then by the decomposition \cite{MarsdenScheurle93}
\begin{equation}
T(Q \times G) \cong TQ \times G \times \mathfrak{g}  \label{trivialization}
\end{equation}
a point in $T^*M$ has coordinates $(q,\dot{q},g,\dot{g})$.
Now consider the \emph{body angular velocity} defined by
\begin{equation}
 \xi = (L_{g^{-1}})_*\,\dot{g},\label{eq:def_angular_velocity} 
\end{equation}
where $(L_{g^{-1}})_*$ denotes the differential of the left translation $L_{g^{-1}}:G \to G, \hspace{3mm} h\mapsto g^{-1}h$, 
at the unit element of $G$.  The commonly used notion of  body angular velocity comes from the fact that in the example where $G$ is the Lie group $SO(3)$ (see Sec.~\ref{sec:Red}) $\xi$ is indeed the angular velocity in a body fixed frame.
Equation~\eqref{eq:def_angular_velocity} is called the \emph{reconstruction equation}
as it can be used to find the full dynamics corresponding to the reduced one.
We will comment on this in more detail in Sec.~\ref{sec:Reconstruction} for the case of three-body systems.

One can see that the body angular velocity is invariant under
the group action:
for  $h \in G$, define the curve
$m(t)=hg\,(t)$, then
\begin{equation}
(L_{m^{-1}})_*\hspace{1mm}\dot{m}=(L_{(g^{-1}h^{-1})})_*\hspace{1mm} ((L_{h})_* \, \dot{g})=\xi.
\end{equation}
As $\xi$ is invariant under the group action,
 the coordinates $(q,\dot{q},\xi)$ 
give a coordinate system on $(TM)/G$ . 

If $L:TM \rightarrow \mathbb{R}$ is a regular Lagrangian function which is invariant under the action of $G$,
then the function $l:TM/G \rightarrow \mathbb{R}$ given by
\begin{equation}
 l(q,\dot{q},\xi):=L(q,\dot{q},g,\dot{g})
\end{equation}
is well-defined.
This is done by passing to the coordinates $(q,\dot{q},\xi)$, and as
the Lagrangian $L$ is invariant, it is possible to put $L$ in 
the form  of the function $l$ in which the $G$ coordinates disappear.
Using $l$ we can define momenta conjugate to $q$ and $\xi$ as
\begin{equation}
 p_q=\frac{\partial l}{\partial \dot{q}}\,, \label{coord1}
\end{equation}
and
\begin{equation}
\eta=\frac{\partial l}{\partial \xi}\,, \label{coord2}
\end{equation}
respectively.
Here $\eta$ is called the \emph{body angular momentum}, and by the chain rule,
\begin{equation}
 \eta=(L_{g})^*\hspace{1mm} p_g \,, \label{rec}
\end{equation}
where $p_g=\partial L/\partial \dot{g}$ is the conjugate momentum of $g\in G$ \cite{Blochetal96}.  
Like body angular velocity the notion body angular momentum again comes from the context of reduction of rotational symmetries.

Let $H:T^*M \rightarrow \mathbb{R}$ be the Hamiltonian obtained from the Legendre transformation of  the Lagrangian $L$ \cite{Arnold78}, i.e., in coordinates
\begin{equation}
 H(q,p_q,g,p_g)=\dot{q}p_q+\dot{g}p_g-L(q,\dot{q},g,\dot{g})\,.
\end{equation} 
As $H$ is $G$ invariant it induces a function $h$ on $T^*M/G$  given by
\begin{equation}
 h(q,p_q,\eta):=H(q,p_q,g,p_g).
\end{equation}
From the construction above one obtains for  $z=(q,p_q,g,p_g)$ and $\zeta \in \mathfrak{g}$,
\begin{equation}
\begin{split}
\langle J(z), \zeta \rangle &= \langle (p_q,p_g), \zeta_M \rangle \\ 
&=\langle (p_q,p_g),(0,(R_g)_* \zeta) \rangle \\
&= \langle p_g , (R_g)_* \zeta \rangle \\
&=\langle (L_{g^{-1}})^* \eta , (R_g)_* \zeta \rangle  \\
&= \langle (Ad_{g^{-1}})^* \eta , \zeta \rangle. \label{eq:momentummap_lastline}
\end{split}
\end{equation}
Here the first equality follows from the definition of the momentum map \eqref{momentum}, 
 the second equality makes use of the decomposition \eqref{trivialization}, the fact that $ \zeta_M$ is tangent to the group orbit which we identify with $G$ and  \eqref{Lie}, the third equality is clear, the fourth equality uses \eqref{rec} and the final equality follows from
 the definition of the coadjoint action in \eqref{eq:def_coadjoint_action} and \eqref{eq:def_adjoint_action}.
We thus obtain
\begin{equation}
 J(p)=(Ad_{g^{-1}})^* \eta
\end{equation}
or equivalently
\begin{equation}
 \eta=(Ad_{g})^*J(p). \label{eq1}
\end{equation}
Let $\mathcal{O}_\mu$ stand for the coadjoint orbit through $J(p)=\mu \in \mathfrak{g}^*$ for some
fixed $\mu \in \mathfrak{g}^*$, i.e.
\begin{equation}
\mathcal{O}_\mu=\{(Ad_{g^{-1}})^* \mu | \hspace{2mm} g\in G \}  \subset \mathfrak{g}^*\,.
\end{equation}
Then Eq.~\eqref{eq1} gives that
\begin{equation}
 \eta \in \mathcal{O}_\mu. \label{eta}
\end{equation}

Now consider the \emph{reduced  space}
\begin{equation}
P_{\mu}:=J^{-1}(\mathcal{O}_{\mu})/G\,.
\end{equation}
By \eqref{eta}, we conclude that if $J(p)=\mu$ for $p=(q,p_q,g,p_g)$ and some fixed $\mu \in \mathfrak{g}^*$, then for 
$ \eta=(Ad_{g})^*J(p)$, we have
$(q,p_q,\eta)\in P_{\mu}$, i.e. $(q,p_q,\eta)$ are coordinates on the reduced space. The reduction to the space $P_{\mu}$ we described  is a coordinate-based form of
the \emph{orbit reduction} of Marle \cite{OrtegaRatiu04}.  For $n$-body systems, the reduction procedure can be interpreted as passing to a \emph{body-fixed frame}  (cf. Sec.~\ref{sec:Red}). 
In fact, the space $P_{\mu}$ is symplectomorphic to the \emph{Marsden-Weinstein reduced space}
$J^{-1}(\mu)/G_\mu$, where $G_\mu$ is the isotropy group of $\mu\in  \mathfrak{g}^*$. For us,  the reduced space $P_{\mu}$  and the coordinates $(q,p_q,\eta)$
form the basis for defining canonical coordinates on the reduced phase space.

%%%%%%%%%%%%%%%%%%%%%%%%%%%%%%%%%%%%%%%%%%%%%%%
\subsection{Canonical coordinates}
\label{canonicalcoordinates}

The coordinates $(q,p_q,g,p_g)$ defined above are clearly canonical, whereas the coordinates $(q,p_q,g,\eta)$ as seen below are not canonical. 
By (\ref{coord1}) and (\ref{coord2})  we get the Poisson bracket equalities 
\begin{equation}
\{q_{\alpha},\eta_{a}\}=\{p_{q_{\alpha}},\eta_{a}\}=0
\end{equation} 
on $T^*M$
for all $\alpha$ and $a$. Now recall the identification $\xi=(L_{g^{-1}})_* \hspace{1mm} \dot{g}$, and let  $(e_1,...,e_l)$ be a basis
of $\mathfrak{g}$, where $l$ is the dimension of $G$. Choosing  $(q_{\alpha},\dot{q}_{\alpha}
,g_a,\xi_a)$ in place  of $(q_{\alpha},\dot{q}_{\alpha},g_a,\dot{g}_a)$ as coordinates on $TM$
requires one to obtain the dynamics in terms of the \emph{anholonomic frame} $(e_1,...,e_l)$. We refer to \cite{LittlejohnReinsch97} for
a concise derivation of the Lagrangian and Hamiltonian in this context.  We briefly give a derivation of the Poisson
brackets in anholonomic frames in Appenix~\ref{anholonomic}. Using \eqref{eq:PoissonAnholonFrame} we get
that
\begin{equation}
 \{\eta_a,\eta_b\}=-\gamma_{ab}^c \, \eta_c,
\end{equation}
 where $\gamma_{ab}^c$ are the structure constants given by $[e_a,e_b]=\gamma_{ab}^c \, e_c$. This is in fact
the same as the $(-)$ Lie-Poisson bracket \cite{Marsden92} on $G$ (see also Appendix~\ref{lie-poisson}).
As the Poisson structure on the coadjoint orbit
 $\mathcal{O}_\mu$ is the reduced one from the Lie-Poisson structure on $G$, the discussion above
suggests that if a canonical coordinate system $(u,v)$ is chosen on $\mathcal{O}_\mu$ such that
\begin{equation}
\eta=\eta(u,v), \label{nu}  
\end{equation}
then the coordinate system $(q,p_q ,u,v)$  becomes a canonical coordinate system on $P_\mu$. 

%What we have obtained so far may be formalized as follows.

%\begin{theorem}
% The coordinate system given by \eqref{coord1}, \eqref{coord2}, \eqref{nu} gives a local canonical coordinate system
%on the reduced symplectic space.
%\end{theorem}

For a more detailed discussion of canonical coordinates on coadjoint orbits, we refer to \cite{Symes80} among others.
In our example of the three body-problem in Sec.~\ref{sec:Red}
the coadjoint orbits are the body-angular momentum spheres.
In this case the canonical coordinates may be chosen as \emph{Deprit coordinates} \cite{Deprit67} on the body-angular momentum sphere
(see \eqref{sphere1} below). 

As for the dynamics, if one passes to the coordinates given in \eqref{nu}, then the reduced Hamiltonian $h_{\mu}:=h |_{P_{\mu}}$ may be written in the form
\begin{equation}
h_{\mu}=h_{\mu}(q,p_q,u,v).
\end{equation}
Then the equations of motion for the
reduced system have the familiar form
\begin{equation}
\begin{split}
\dot{z}=\{z_,h_\mu\}
\end{split}
\end{equation}
for $z=(q,p_q,u,v)$.

%%%%%%%%%%%%%%%%%%%%%%
\subsection{Simple mechanical systems}
\label{sec:simplemechsyst}

The construction of canonical coordinates on the reduced phase space and the reduction of  the equations of motion can be described more explicitly in the case of simple mechanical systems. 
For a simple mechanical system, the Lagrangian is of the form
\begin{equation}
 L(s,\dot{s})=\frac{1}{2}\dot{s}^T \, k \, \dot{s}-V(q),
\end{equation}
where $s$ is a coordinate system on $M$, $k$ is a Riemannian metric on $M$ which is invariant under the action of $G$ on $M$, and $V$ is a potential function which is also invariant and hence, if 
$(q,g)$ is a coordinate system on the local trivialization
  \begin{equation}
 M\cong Q \times G \label{trivial}
\end{equation}
then $V$
depends only on the shape coordinates $q$.
Let  $\sigma :Q \rightarrow M$ be a local section of the fibre bundle, i.e.
$\sigma$ is a right inverse of the projection $\pi :M \rightarrow Q$.
If a point $q\in Q$ is given, then a point $s\in \pi ^{-1}(q)$ is of the form $g \sigma (q)$ with some $g\in G$. Set $r=\sigma(q)$, then
\begin{equation}
 \dot{r}=\frac{\partial r}{\partial q} \hspace{1mm}\dot{q}
\end{equation}
by the chain rule.

Let the \emph{body velocities} be  defined as
\begin{equation}
 v=(L_{g^{-1}})_*\dot{s}\,. \label{eq:def_body_velocities}
\end{equation}
By using the Leibniz rule \cite{AbrahamMarsden78}
and the definition of the fundamental vector field
one obtains
\begin{equation}
v=\xi_M(r)+\dot{r}, \label{eq:decomp_v}
\end{equation}
where $\xi_M$ is the fundamental vector field corresponding to $\xi=(L_{g^{-1}})_* \hspace{1mm} \dot{g}\in \mathfrak{g}$. 
The kinetic energy thus becomes
\begin{equation}
K=\frac{1}{2} \, v^T \, k \, v=\frac{1}{2} \, \xi_M(r) ^T \, k \, \xi_M(r)+\,
\xi_M(r) ^T \, k \, \dot{r}
+\frac{1}{2} \, \dot{r}^T \, k \, \dot{r}. 
\end{equation}
For $\xi,\eta\in \mathfrak{g}$, set
\begin{equation}
\xi^T \, \mathbb{I} \, \eta = \xi_M(r)^T \, k \, \eta_M(r)\,. 
\end{equation}
Then    $\mathbb{I}$ is
a left-invariant inner product on $\mathfrak{g}$ \cite{Marsden92}. For $G=SO(3)$ (respectively $ \mathfrak{g}=so(3)$), $\mathbb{I}$ is the moment of inertia tensor (see Sec.~\ref{sec:Red}).
In order to decouple the kinetic energy in group and shape terms, 
the so called \emph{mechanical
connection} is introduced \cite{Marsden92}. The mechanical connection
$A:TM \rightarrow \mathfrak{g}$ is defined by
\begin{equation}
A(s,\dot{s})=\mathbb{I}^{-1} \mathbf{L}(s,\dot{s}),
\end{equation} 
where $\mathbf{L}$ is the tangent momentum map given in \eqref{tan_mom_map},
and $\mathbb{I}^{-1}:\mathfrak{g}^* \rightarrow \mathfrak{g}$
is the linear map associated with the inner product $\mathbb{I}$.
At any point $s\in M$, the tangent space to $T_sM$ may be decomposed into
\begin{equation}
 T_sM=V_s+H_s,
\end{equation}
where $V_s$ is the tangent space to the orbit $Gs$, and $H_s$ is the space which is orthogonal
to $V_s$ with respect to the metric $k$.
A tangent vector $w \in T_sM$ may be written in this decomposition as
\begin{equation}
 w=\text{ver}_s w + \text{hor}_s w.
\end{equation}
It turns out that \cite{Marsden92}
\begin{equation}
 \text{ver}_s w = [A(s,w)]_M(s)
\end{equation}
and
\begin{equation}
 \mathbf{L}(s,\text{hor}_s w)=0.
\end{equation}
With respect to this decomposition $\dot{r}$
may be written in the form
\begin{equation}
 \dot{r}=\text{ver}\hspace{1mm} \dot{r}+(\dot{r}-\text{ver}\hspace{1mm} \dot{r}).
\end{equation}
If a new metric $d$ is introduced by 
\begin{equation}
 \dot{q}^T \, d \, \dot{q}:=(\dot{r}-\text{ver}\hspace{1mm} \dot{r})^T \,
 k \, (\dot{r}-\text{ver}\hspace{1mm} \dot{r}), \label{eq:def_horizontal_metric}
\end{equation}
which is called the \emph{horizontal metric},
and if we define the map $A_Q:TQ\to  \mathfrak{g}$ by
\begin{equation}
 A_Q \hspace{1mm} \dot{q}:=A \hspace{1mm}\dot{r},
\end{equation}
then
after rearranging terms the kinetic energy assumes the form
\begin{equation}
 K=\frac{1}{2} \, (\xi+A_Q\hspace{1mm} \dot{q})^T \,
 \mathbb{I} \, (\xi+A_Q\hspace{1mm} \dot{q})+\frac{1}{2} \, \dot{q}^T\,d\,\dot{q} .
\end{equation}
This compact form of the kinetic energy reflects the decomposition of it into vertical and horizontal
energies. Finally the Lagrangian in coordinates $(q,\dot{q},\xi)$ is given by

\begin{equation}
 l(q,\dot{q},\xi)=\frac{1}{2} \, (\xi+A_Q\hspace{1mm} \dot{q})^T \,
  \mathbb{I} \, (\xi+A_Q\hspace{1mm} \dot{q})+\frac{1}{2} \, \dot{q}^T\,d\,\dot{q} - V(q).
\end{equation}
If the \emph{body angular momentum} is defined by 
\begin{equation}
 \mathbf{J}=(L_{g^{-1}})_* \mathbf{L}, \label{eq:def_angular_momentum_simple}
\end{equation}
then it is seen that
\begin{equation}
 \mathbf{J}=\mathbb{I} \, (\xi + A_Q \hspace{1mm} \dot{q}).
\end{equation}
 On the other hand, the conjugate momenta of $q$ and $\xi$  are obtained as
\begin{equation}
p_q =\frac{\partial l}{\partial \dot{q}}=d \, \dot{q}+A_Q^T \, \mathbf{J},
\end{equation}
and
\begin{equation}
 \eta=\frac{\partial l}{\partial \xi}=\mathbb{I} \,( \xi + A_Q \hspace{1mm} \dot{q}),
\end{equation}
respectively. Note here that $\eta=\mathbf{J}$. Finally
the Hamiltonian can be written as
\begin{equation}
 h(q,p_q,\mathbf{J})=\frac{1}{2} \, \mathbf{J}^T \, \mathbb{I}^{-1} \, \mathbf{J} +
 (p_q- A_Q^T \, \mathbf{J} )^T \,
 d^{-1} \, (p_q- A_Q^T \, \mathbf{J})
+ V(q), \label{simplered}
\end{equation}
where $d^{-1}$ denotes the metric on $T^*M$ corresponding to $d$. 

%%%%%%%%%%%%%%%%%%%%%%%%%%%%%%%%%%%%%%%%%%%%%%%%%%
\subsection{Special cases}
\label{sec:examples}

\textbf{(1) $G$ is Abelian.}  Consider the case where the Lie group is Abelian, e.g. a torus group. The reduction strongly simplifies in this case. Since
the coadjoint action is trivial  the well-known identification \cite{OrtegaRatiu04}
\begin{equation}
J^{-1}(\mathcal{O}_{\mu})/G = T^{*}Q 
\end{equation}
is obtained. This shows that one can take the coordinates $(q^\alpha,p_q ^\alpha)$ as the canonical coordinates on the reduced space,
which is symplectomorphic to $T^{*}Q$. %, so the normal form procedure is quite obvious.
An example of this situation is given by  the translational motions of an $n$-body system \cite{LittlejohnReinsch97}.
We will explicitly illustrate the Abelian case for the example of a double spherical pendulum in Sec.~\ref{sec:Pen}.
A detailed analysis of the Abelian case in the Lagrangian setting
can be found in \cite{MarsdenScheurle93}.

\textbf{(2) Vanishing angular momentum.} With the notation above, if $\eta \equiv 0$, then the coadjoint orbit is trivial
 as in the first special case. Then the reduced space is symplectomorphic to $T^{*}Q$. A well studied example of a system with vanishing angular momentum is 
the so called \emph{falling cat problem} \cite{Montgomery93}.

\textbf{(3) Generalized rigid bodies.} Suppose that $M=G$. Then $Q$ is just a point and 
\begin{equation}
J^{-1}(\mathcal{O}_{\mu})/G = \mathcal{O}_{\mu}.
\end{equation}
This  occurs,  e.g., for a rigid body, where the configuration space $M$ is the rotation group $SO(3)$
and the reduced space is the body-angular momentum sphere. 

%%%%%%%%%%%%%%%%%%%%%%%%%%%%%%%%%%%%%%%%%%%%%%%%%%%%%%%
\subsection{Relative Equilibria}
\label{sec:relative}
A point on $T^*M$ is called a \emph{relative equilibrium point} if its projection into
the reduced space is a critical point of the reduced Hamiltonian. So we are interested in the
equilibria of the function $h_\mu :P_{\mu}\rightarrow \mathbb{R}$
for some fixed $\mu \in \mathfrak{g}^*$.
We will give some criteria for relative equilibria for simple mechanical systems.

From Eq.~\eqref{simplered} the equations of motion are obtained to be
\begin{equation}
\begin{split}
 \dot{\mathbf{J}}=&-\mbox{ad}^*_{\partial h /\partial \mathbf{J}} \, \mathbf{J} \\
 \dot{q}=&\frac{\partial h}{\partial p_q}=d^{-1}\,(p_q-A_Q^T \, \mathbf{J}), \\
 \dot{p_q}=&-\frac{\partial h}{\partial q}=-\frac{1}{2}\, \frac{\partial}{\partial q}\big(
  (p_q- A_Q^T \, \mathbf{J})^T \, d^{-1} \, (p_q-A_Q^T \, \mathbf{J})\big)
  +\frac{\partial V_{\text{eff}}}{\partial q},
\end{split}
\end{equation}
where 
\begin{equation}
\frac{\partial h}{\partial \mathbf{J}}=\mathbb{I}^{-1} \, \mathbf{J}-A_Q \, d^{-1}
(p_q-A_Q^T \, \mathbf{J}) 
\end{equation}
and
\begin{equation}
V_{\mbox{eff}}=\frac{1}{2} \, \mathbf{J}^T \, \mathbb{I}^{-1} \,\mathbf{J}
+ V(q)
\end{equation}
 is the \emph{effective potential}. We note that for Abelian actions, the effective potential agrees with 
the so called \emph{amended potential} \cite{Marsden92}.
The 
equations above can be deduced from Appendix~\ref{anholonomic} and Appendix~\ref{lie-poisson}.
(See also \cite{LittlejohnReinsch97,KozinRobertsTennyson00}
and \cite{Cendraetal03}.)
Then the conditions for having a relative equilibrium are 
\begin{equation}
\begin{split}
p_q=& A_Q^T \, \mathbf{J}, \\
\mbox{ad}^*_{ \mathbb{I}^{-1} \, \mathbf{J}} \, \mathbf{J}  =&0\,,\\
 \frac{\partial}{\partial q} V_{\text{eff}} =&0\,.
 \end{split} \label{eq:cond_rel_equil_gen}
 \end{equation}

In the examples of the three-body problem and the double spherical pendulum  in Sections~\ref{sec:Red} and \ref{sec:Pen}
the conditions will be given in a more explicit form.

%%%%%%%%%%%%%%%%%%%%%%%%%%%%%%%%%%%%%%%%%%%%%
\subsection{Poincar{\'e}-Birkhoff normal form around a relative equilibrium point}
\label{sec:normal}

The Poincar{\'e}-Birkhoff normal form is a main tool in dynamical systems theory.  It allows one (if certain conditions are satisfied) to study the dynamics of a nonlinear system in the neighborhood of an equilibrium point by approximating it by a `simpler' system. This has many applications, e.g., in the study of bifurcations and the computations of center manifolds \cite{Deprit69,AKN88,DragtFinn76,Murdock03}.  
The simpler system is constructed order by order of the Taylor expansion of the original system at the equilibrium point by a suitable choice of coordinates at each order. 
For Hamiltonian systems,  the coordinate transformation are sought to be symplectic.  As the dynamics  (i.e. the vector field)  is generated by a Hamilton function the simplification can be described completely in terms of a simplification of the Hamilton function.
There are well established  algorithms which can be implemented on a computer and which allow one to compute normal forms to any desired order (see Appendix~\ref{algorithm}).
As the starting point for these algorithms is a Hamiltonian system with canonical coordinates on the linear symplectic space  $\mathbb{R}^f\times \mathbb{R}^f$ where $f$ denotes the number of degrees of freedom it is crucial for the application of these algorithms to relative equilibria of symmetry reduced Hamiltonian systems to explicitly construct canonical coordinates on the reduced space as described in the subsections above.  

We will in this paper restrict ourselves to Poincar{\'e}-Birkhoff normal forms at  equilibrium points 
where the eigenvalues associated with  the linearized Hamiltonian vector field $\mathfrak{J} D^2 H$ are purely imaginary. Here $\mathfrak{J}$ denotes the standard symplectic matrix and $D^2 H$ is the Hessian of the Hamiltonian $H$.
We will denote the eigenvalues by $\pm \ui \,\omega_k$, $k=1,\ldots,f$.
 Assuming that the eigenvalues are independent over the field of rational numbers
(i.e. in the absence of \emph{resonances}), 
the Poincar{\'e}-Birkhoff normal form yields a symplectic transformation to new (\emph{normal form}) coordinates such that 
the transformed Hamiltonian function truncated at order $n_0$ of its Taylor expansion assumes the  form 
\begin{equation} 
H_{\text{NF}}(I_1,\ldots,I_\dof) = \sum_{k=1}^f \omega_k I_k + \text{h.o.t.} \,,
\end{equation}
where $I_k$, $k=1,\ldots,\dof$, are constants of motions which (when expressed in terms of the normal form coordinates)
have the form
\begin{equation}
 I_k=p_k^2+q_k^2, \hspace{3mm} k=1,\ldots,f,
\end{equation}
and $H_{\text{NF}}$ is a polynomial 
of order $n_0/2$ in  $I_k$, $k=1,\ldots,\dof$ and hence of order $n_0$ in $p$ and $q$ (note that only even orders $n_0$ of a normal form make sense).
The algorithm to compute this transformation is sketched in Appendix~\ref{algorithm}. We will apply it to the examples of relative equilibria of a three-body system and a double spherical pendulum in Sec.~\ref{sec:Normal} and Sec.~\ref{stretched}, respectively.

% end rem

%%%%%%%%%%%%%%%%%%%%%%
\section{Three-body systems}
\label{sec:Red}
In this section we review the reduction of a three-body system for which we then write the reduced Hamiltonian in canonical coordinates following Sec.~\ref{sec:2}. We also comment on the reconstruction of the full (unreduced) dynamics. As an example we discuss a  triatomic molecule with  a Morse-type potential for which we compute the Poincar{\'e}-Birkhoff 
normal form about an equilibrium point given by an equilateral triangle configuration.

%%%%%%%%%%%%%%%%%%%%%
\subsection{Reduced equations of motion}

\begin{figure}
\begin{center}
\includegraphics[angle=0,width=7cm]{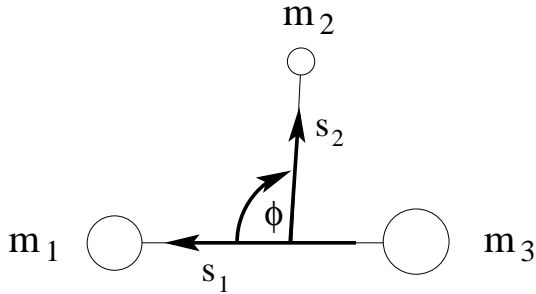}
\end{center}
\caption{\label{fig:conf}
Definition of the Jacobi vectors $\mathbf{s}_{1}$ and $\mathbf{s}_{2}$ and the corresponding angle $\phi$.
}
\end{figure}

Consider a  system of three bodies with masses $m_1,m_2,m_3$ and position vectors $\mathbf{x}_{1},\mathbf{x}_{1},\mathbf{x}_{1} \in \mathbb{R}^3$,
respectively, without external forces acting on the three bodies. The symmetry of overall translations can be reduced by introducing
mass-weighted Jacobi vectors which are defined according to  
\begin{eqnarray*}
\mathbf{s}_1&=&\sqrt{\mu_1}(\mathbf{x}_{1}-\mathbf{x}_{3}),\\
\mathbf{s}_2&=&\sqrt{\mu_2} 
(\mathbf{x}_{2}-\frac{m_{1}\mathbf{x}_{1}+m_{3}\mathbf{x}_{3}}{m_{1}+m_{3}}),
\end{eqnarray*}%
where
\begin{equation}
 \mu_1=\frac{m_{1}m_{3}}{m_{1}+m_{3}}, \ \ \ \mu_2=\frac{m_{2}(m_{1}+m_{3})}{m_{1}+m_{2}+m_{3}}
\end{equation}
are  reduced masses (see Figure \ref{fig:conf}). 

Excluding collinear (and hence also collisional) configurations 
we obtain the six-dimensional translation-reduced configuration space 
\begin{equation}
M=\left\{ s=(\mathbf{s}_1,\mathbf{s}_2):\ \lambda \mathbf{s}_1+\mu \mathbf{s}_2\neq
0 \mbox{ for all } (\lambda ,\mu) \in \mathbb{R}^2\backslash \{0\}\right\} \subset \mathbb{R}^{3}%
\mathbb{\times R}^{3}.
\end{equation}%
Proper rotations $g\in SO(3)$ act on $M$ in the natural way
\begin{equation}
g(\mathbf{s}_1,\mathbf{s}_2)=(g\mathbf{s}_1,g\mathbf{s}_2).
\end{equation}%
On $M$ this action is free and it thus follows that the shape space
\begin{equation}
 Q:=M/SO(3)
\end{equation}
has a manifold structure which turns out to be diffeomorphic to
$\mathbb{R}^3_{+}=\{(x,y,z)\in \mathbb{R}^3 \,:\,  z>0\}$ \cite{Iwai87c}.
The canonical projection $M\rightarrow Q$ defines a principal bundle with structure group 
$SO(3)$. This principal bundle is trivial \cite{Iwai87c}
and has holonomy group $SO(2)$ yielding a geometric reduction 
\cite{CiftciWaalkens10}.

The Lie algebra $\mathfrak{g}=\mathfrak{so}(3)$ of $G=SO(3)$
can be identified with $\mathbb{R}^3$ where the Lie algebra structure
becomes the vector product `$\times$'. By using the
bi-invariant inner product on $\mathfrak{g}$, or
equivalently  
the dot product on $\mathbb{R}^3$,
one can identify $\mathfrak{g}^*$ also with $\mathbb{R}^3$. 
With these identifications
the fundamental vector field
corresponding to $\zeta \in \mathfrak{g}$ at $s=(\mathbf{s}_1,\mathbf{s}_2)$ is
\begin{equation}
\zeta(s)=(\zeta \times \mathbf{s}_1,\zeta \times \mathbf{s}_2)\,.
\end{equation}
The momentum mapping  $\mathbf{L}:TM\to \mathfrak{g}^*$,  following \eqref{tan_mom_map}, is given by
\begin{equation}
\mathbf{L}= \mathbf{s}_1 \times \mathbf{\dot{s}}_1+\mathbf{s}_2 \times \mathbf{\dot{s}}_2. 
\end{equation}

After choosing a body-fixed frame one can obtain the corresponding body-fixed Jacobi vectors by
\begin{equation}
 \mathbf{s}_i=g\, \mathbf{r}_i, \hspace{3mm} i=1,2\,,
\end{equation}
where $g\in SO(3)$ is the matrix relating the body-fixed frame and the space-fixed frame.  As $g$ depends on three coordinates, e.g. Euler angles, 
there are three shape space coordinates $q_\alpha$, $\alpha=1,2,3$,  remaining to parametrize the two vectors $\mathbf{r}_1$ and $\mathbf{r}_2$.  

The kinetic energy is given by
\begin{equation}
K=\frac{1}{2} \sum^{3}_{i=1} m_i \, \mathbf{\dot{x}}^{2}_i\,=\frac{1}{2} \sum^{2}_{i=1} \mathbf{\dot{s}}^{2}_i\,. 
\end{equation}
The corresponding metric $k$ thus is Euclidean.
Defining body velocities according to
\begin{equation}
 \mathbf{v}_i=g^T \, \mathbf{\dot{s}}_i, \hspace{3mm} i=1,2,
\end{equation}
(cf. \eqref{eq:def_body_velocities}), 
and using the shape coordinates and their time derivatives one can rewrite the body velocities as
\begin{eqnarray*}
\mathbf{v}_i&=&g^T \, (\dot{g} \, \mathbf{r}_i  
+\sum^{3}_{\alpha=1} g \, \frac{\partial \mathbf{r}_i}{\partial q_{\alpha}}\dot{q}_\alpha) \\
&=&g^T \, \dot{g} \, \mathbf{r}_i +\sum^{3}_{\alpha=1} \frac{\partial \mathbf{r}_i}{\partial q_{\alpha}}\, \dot{q}_\alpha.
\end{eqnarray*}
The body angular velocity $\mathbf{\xi}$ (see \eqref{eq:def_angular_velocity}) is the vector
in $\mathbf{R}^3$ corresponding to $\mathbf{\Xi} \in \mathfrak{g}$
given by
\begin{equation}
 \mathbf{\Xi}=g^T \, \dot{g}\,, \label{Omega}
\end{equation}
which is the reconstruction equation \eqref{eq:def_angular_velocity}.
Then one has
\begin{equation}
\mathbf{v}_i=\mathbf{\xi} \times \mathbf{r}_i+\sum^{3}_{\alpha=1} \frac{\partial \mathbf{r}_i}{\partial q_\alpha}\dot{q}_\alpha \label{velocity}
\end{equation}
which corresponds to the general expression \eqref{eq:decomp_v}.
Since the moment of inertia tensor $\mathbb{I}$ is given by 
\begin{equation}
 \mathbb{I} \, \mathbf{u}=\mathbf{r}_1 \times (\mathbf{u} \times \mathbf{r}_1)+\mathbf{r}_2 \times (\mathbf{u} \times \mathbf{r}_2),
\end{equation}
for $\mathbf{u}\in \mathbb{R}^3$, the mechanical connection
$A_Q=\left[\begin{array}{ccc}
           \mathbf{A}_1\\
\mathbf{A}_2\\
\mathbf{A}_3
           \end{array}
\right]$ 
is obtained to be
\begin{equation}
 \mathbf{A_{\alpha}}=\mathbb{I}^{-1}(\mathbf{r}_1 \times \frac{\partial \mathbf{r}_1}{\partial q_\alpha})+
\mathbb{I}^{-1}(\mathbf{r}_2 \times \frac{\partial \mathbf{r}_2}{\partial q_\alpha}).
\end{equation}
Then the kinetic energy  becomes
\begin{equation}
K=\frac{1}{2} \, \mathbf{\xi}^T \, \mathbb{I} \, \mathbf{\xi}+\sum^{3}_{\alpha=1} (\mathbf{\xi}^T \, \mathbb{I} \, \mathbf{A_\alpha}) \,\dot{q}_\alpha
+\frac{1}{2} \, \sum^{3}_{\alpha,\beta=1}h_{\alpha \beta}\, \dot{q}_\alpha \, \dot{q}_\beta,   
\end{equation}
where 
\begin{equation}
h_{\alpha \beta}=\sum^{2}_{i=1} \frac{\partial \mathbf{r}_i}{\partial q_\alpha}^T \, \frac{\partial \mathbf{r}_i}{\partial q_\beta}.
\end{equation}
Using that the horizontal metric is 
\begin{equation}
 d_{\alpha \beta}=h_{\alpha \beta}-\mathbf{A}_{\alpha}^T \, \mathbb{I} \, \mathbf{A_\beta}
\end{equation}
(see \eqref{eq:def_horizontal_metric})
allows one to write  the kinetic energy in the  compact form
\begin{equation}
 K=\frac{1}{2}\sum^{3}_{\alpha,\beta=1}(\mathbf{\xi}+\mathbf{A_\alpha} \dot{q}_\alpha)^T \, \mathbb{I} \, (\mathbf{\xi}+
\mathbf{A_\beta} \dot{q}_\beta)+\frac{1}{2} \,d_{\alpha\beta}\dot{q}_\alpha \dot{q}_\beta.
\end{equation}

Following \eqref{eq:def_angular_momentum_simple} the body angular momentum is given by
\begin{equation}
 \mathbf{J}=g^T \, \mathbf{L}= \mathbf{r}_1 \times \mathbf{v}_1+\mathbf{r}_2 \times \mathbf{v}_2.
\end{equation}
Then by Eq.~\eqref{velocity} and one has
\begin{equation}
 \mathbf{J}= \mathbb{I} \cdot (\mathbf{\xi +\sum^{3}_{\alpha=1} \mathbf{A}_\alpha} \,\dot{q}_\alpha). \label{angular}
\end{equation}
The conjugate momenta are given by
\begin{equation}
 \eta=\frac{\partial K}{\partial \xi}=\mathbb{I} \cdot (\mathbf{\xi} +\sum^{3}_{\alpha=1} \mathbf{A}_\alpha \, \dot{q}_\alpha)=\mathbf{J},
\end{equation}
and
\begin{equation}
 p_\alpha =\frac{\partial K}{\partial q_\alpha}=\sum^{3}_{\beta=1}d_{\alpha \beta} \dot{q}_\beta + \mathbf{J}^T \, \mathbf{A}_\alpha. \label{three_moment}
\end{equation}
Thus the Hamiltonian takes the form
\begin{equation}
h=\frac{1}{2}\,\mathbf{J}^T \, \mathbb{I}^{-1} \, \mathbf{J}+
\frac{1}{2} \sum^{3}_{\alpha,\beta=1} d^{\alpha \beta} (p_{\alpha}-\mathbf{J}^T \, \mathbf{A}_\alpha)(p_\beta-\mathbf{J}^T \, \mathbf{A}_\beta)+V, \label{Ham}
\end{equation}
where $V=V(q_1,q_2,q_3)$ is the potential.

\begin{figure}
\begin{center}
\includegraphics[angle=0,width=7cm]{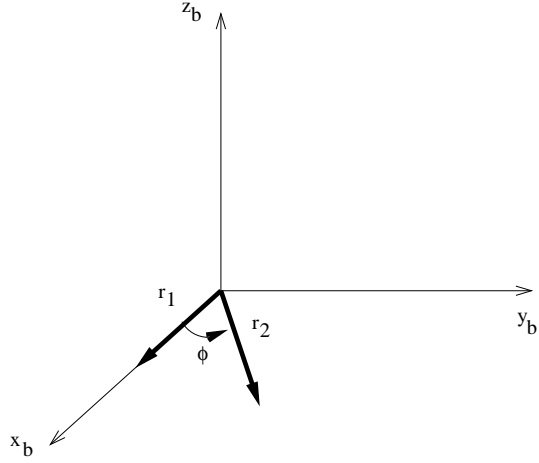}
\end{center}
\caption{\label{fig:xxy}
Definition of a body-fixed frame according to the $xxy$-gauge.
}
\end{figure}

Let us now make give explicit expression by introducing coordinates. 
As the shape coordinates $(q_1,q_2,q_3)$ we choose  \emph{Jacobi coordinates} ($r_1,r_2,\phi$) which are defined as \cite{LittlejohnReinsch97}
\begin{equation}
r_1=\sqrt{\mathbf{r_1} \cdot \mathbf{r_1}}, \quad r_2=\sqrt{\mathbf{r_2} \cdot \mathbf{r_2}}, \quad \phi =\cos ^{-1} (\mathbf{r_1} \cdot \mathbf{r_2} / (r_1r_2)), \quad 0 \leq \phi \leq \pi, \label{coor}
\end{equation}
(see Fig.~\ref{fig:conf}). Choosing then the axes $x_b,y_b,z_b$ of a body-fixed frame according to the so called $xxy$-gauge\footnote{We note that the choice of a body-fixed frame corresponds to the choice of the local 
section $\sigma$ of the fibre bundle $M\to Q$ in Sec.~\ref{sec:simplemechsyst}. The gauge theoretical interpretation of this choice is studied in great detail in \cite{LittlejohnReinsch97}.} 
shown in  Fig.~\ref{fig:xxy},
one obtains for the moment of inertia tensor, metric and mechanical connection \cite{LittlejohnReinsch97}
\begin{equation}
\mathbb{I}=\left[
\begin{array}{ccc} 
r_{2}^{2} \sin^2 \phi & -r_{2}^{2} \sin \phi \cos \phi & 0 \\ 
-r_{2}^{2} \sin \phi \cos \phi & r_{1}^{2}+r_{2}^{2} \cos^2  \phi & 0 \\ 
0 & 0 & r_{1}^{2}+r_{2}^{2} 
\end{array}
\right]\,, \label{eq:inertia_tensor_Jacobi}
\end{equation}
\begin{equation}
 \left[ d_{\mu\nu}\right] =\left[ 
\begin{array}{ccc}
1 & 0 & 0 \\ 
0 & 1 & 0 \\ 
0 & 0 & \frac{r_{1}^{2}r_{2}^{2}}{r_{1}^{2}+r_{2}^{2}}%
\end{array}%
\right]\,,
\end{equation}
and
\begin{equation}
 \mathbf{A}_{r_1 }=\mathbf{A}_{r_2}=(0,0,0),\ \ \ \mathbf{A}_{\phi }=(0,0,\frac{r_{2}^{2}}{r_{1}^{2}+r_{2}^{2}} )\,, \label{eq:mech_connection_Jacobi}
\end{equation}
respectively.

\begin{figure}
\begin{center}
\raisebox{8cm}{}\includegraphics[angle=0,width=7cm]{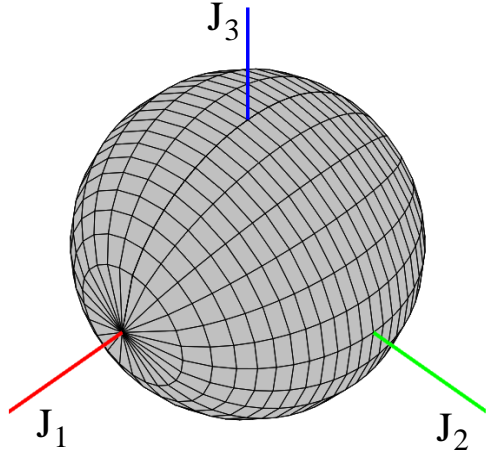}
\end{center}
\caption{\label{fig:Jsphere}
Coordinate lines on the angular momentum sphere of $(u,v)$ defined according to  \eqref{sphere1}.
}
\end{figure}

Putting the results above together the Hamiltonian in terms of Jacobi coordinates becomes
\begin{eqnarray*}
 h(r_1,r_2,\phi,p_1,p_2,p_3,\mathbf{J})&=&\frac{1}{2} \{ \frac{r_{1}^{2}+r_{2}^{2} \cos^2  \phi}{r_{1}^{2}r_{2}^{2} \sin^2  \phi}J_{1}^{2}+
\frac{2 \cos \phi}{r_{1}^{2} \sin \phi}J_{1}J_{2}+\frac{1}{r_{1}^{2}}J_{2}^{2}+\frac{1}{r_{1}^{2}+r_{2}^{2}}J_{3}^{2} \\
&+&p_{1}^{2} +p_{2}^{2}+\frac{r_1^2+r_2^2}{r_1^2 r_2^2}(p_{3}-\frac{r_{2}^{2}}{r_{1}^{2}+r_{2}^{2}}J_{3})^{2} \}+V(r_1,r_2,\phi)\,,
\end{eqnarray*}
where 
$
 \mathbf{J}=(J_1,J_2,J_3)\,.
$
Here $\| \mathbf{J} \|$ is conserved so the coadjoint orbit is the body angular momentum sphere $S^2(\| \mathbf{J} \|)$.
One choice of canonical coordinates on $S^2(\| \mathbf{J} \|)$ are the so called Deprit coordinates which are defined as \cite{Deprit67}
\begin{equation}
 (J_{1},J_2,J_3) = (v,\sqrt{r^{2}-v^{2}}\sin u,\sqrt{r^{2}-v^{2}}\cos u) \label{sphere1}\\
\end{equation}
are chosen, where $r=\| \mathbf{J} \|$ (see Fig.~\ref{fig:Jsphere}), then the reduced Hamiltonian $h|_{S^2(r)}$ becomes
\begin{eqnarray*}
 h_r(r_1,r_2,\phi,p_1,p_2,p_3,u,v)&=&\frac{1}{2} \{ \frac{r_{1}^{2}+r_{2}^{2} \cos^2  \phi}{r_{1}^{2}r_{2}^{2} \sin^2  \phi}v^{2}
+\frac{2 \cos \phi}{r_{1}^{2} \sin \phi}v\sqrt{r^2-v^2}\sin u\\
&+&\frac{1}{r_{1}^{2}}(r^2-v^2)\sin^2 u+\frac{1}{r_{1}^{2} 
+r_{2}^{2}}(r^2-v^2)\cos^2 u+p_{1}^{2}+p_{2}^{2}\\
&+&\frac{r_1^2+r_2^2}{r_1^2 r_2^2}(p_{3}-\frac{r_{2}^{2}}{r_{1}^{2}+r_{2}^{2}}\sqrt{r^2-v^2}\cos u)^{2} \}+V.
\end{eqnarray*}

%%%%%%%%%%%%%%%%%%%%%%%%
\subsection{Lagrangian equilateral triangle configurations}
\label{Lagrangian}

We now consider a three-body system with a  Morse-type potential given by
\begin{equation}
 V=\sum^{3}_{1 \le i < j \le 3} \exp(-2 (r_{ij}-d_0))-2 \exp(-(r_{ij}-d_0))\,, \label{eq:def_Morse_potential}
\end{equation}
where $r_{ij}$ is the distance between the $i$th and the $j$th particle. 
The parameter $d_0$ determines the side length of the equilateral triangle at which the potential has a minimum.  
In Jacobi coordinates one has
\begin{eqnarray}
r_{13}&=&\frac{r_1}{\sqrt{\mu_1}}, \\
r_{23}&=&\sqrt{\frac{\mu_1 r_1^2}{m_3^2}+\frac{r_2^2}{\mu_2}+\frac{2 \sqrt{\mu_1}r_1r_2 \cos \phi}{m_3\sqrt{\mu_2}}}, \\
r_{12}&=&\sqrt{\frac{\mu_1 r_1^2}{m_1^2}+\frac{r_2^2}{\mu_2}-\frac{2 \sqrt{\mu_1}r_1r_2 \cos \phi}{m_1\sqrt{\mu_2}}}.
\end{eqnarray}
A Lagrangian equilateral triangle relative equilibrium is a planar motion where the shape is a constant equilateral triangle. The angular momentum is orthogonal to the plane of the motion. The  
mass-weighted Jacobi vectors are then of the form
\begin{equation}
\begin{split}
\mathbf{r}_{1}^{e}&=\sqrt{\mu_1}(b,0,0), \\
\mathbf{r}_{2}^{e}&=\sqrt{\mu_2}\left(\frac{b}{2}\left(\frac{m_3-m_1}{m_1+m_3}\right),\pm \frac{\sqrt{3}}{2} b,0 \right)\,, \label{eq:Jacobi_Lagrange_triangle}
\end{split}
\end{equation}
where the parameter $b$ is determined by the magnitude of the angular momentum $r$ or conversely, choosing a value for $b$ determines $r$.
The corresponding
Jacobi coordinates  $(r_1^e,r_2^e,\phi^e)$ are easily computed using \eqref{coor}.

Now we find the values of the other coordinates and the parameter $r$ at the equilibria specified by $b$.
Following \eqref{Ham} a relative equilibrium satisfies \cite{KozinRobertsTennyson00}
\begin{eqnarray}\label{Rel}
\mathbf{J}\times \left( \mathbb{I}^{-1}\cdot \mathbf{J}\right)  &=&0, \label{eigen} \\
p_{\alpha } &=&\mathbf{J}\cdot \mathbf{A}_{\alpha }, \label{eq:momenta_rel_equilibria_Jacobi}  \\ 
\frac{\partial }{\partial q_{\alpha }}(\frac12 \mathbf{J}^T \mathbb{I}^{-1}\, \mathbf{J}+V) &=&0. \label{relative}
\end{eqnarray}
By Eq.~\eqref{eigen} $\mathbf{J}$ is an eigenvector of $\mathbb{I}^{-1}$ at a relative equilibrium point. 
For a Lagrangian equilateral triangle relative equilibrium we know that in the $xxy$-gauge $\mathbf{J}$ is pointing in the $z$-direction of the body frame. Hence $\mathbf{J}=(0,0,r)$. From \eqref{sphere1} we find the corresponding canonical coordinates  $(u^e,v^e)=(0,0)$.
Inserting $\mathbf{J}=(0,0,r)$  and using  
the block structure of the inertia tensor in \eqref{eq:inertia_tensor_Jacobi} Eq.~ \eqref{relative}  reduces to 
\begin{equation}
\frac{\partial }{\partial q_{\alpha }}( \frac12 \frac{r^2}{r_1^2+r_2^2}+V)=0 \label{relative_eq}\,.
\end{equation}
We use this equation to find the magnitude of the angular momentum $r$ for a  Lagrangian equilateral triangle $(r_1^e,r_2^e,\phi^e)$ specified by a given parameter $b$ in \eqref{eq:Jacobi_Lagrange_triangle}.
Finally, inserting \eqref{eq:mech_connection_Jacobi} in \eqref{eq:momenta_rel_equilibria_Jacobi}
the conjugate momenta are obtained to be
\begin{equation}
 p_1^e=p_2^e=0, \hspace{3mm} p_3^e=\frac{(r_2^e)^2}{(r_1^e)^2 + (r_2^e)^2}\,r\,.
\end{equation}

\begin{figure}
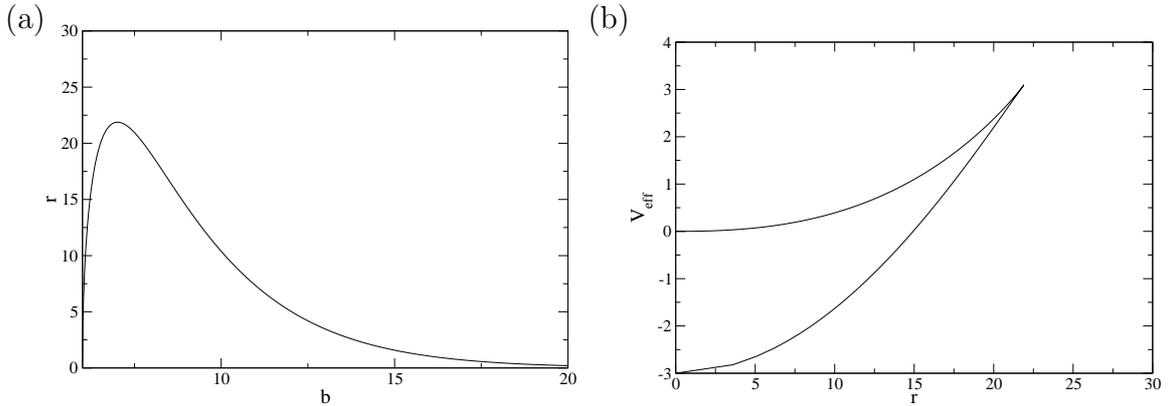

\begin{center}
\raisebox{5cm}{(a)}\includegraphics[angle=0,width=7cm]{three_body_r_versus_b}
\raisebox{5cm}{(b)}\includegraphics[angle=0,width=7cm]{three_body_V_versus_r}
\end{center}
\caption{\label{fig:three_body_angularmomentum_energy}
(a) The magnitude $r$ of the angular momentum  as a function of the size of the Lagrangian equilateral triangle parametrized by $b$ (see \eqref{eq:Jacobi_Lagrange_triangle}).
(b) The energy as given by the effective potential $V_\text{eff}=\frac12 \frac{r^2}{r_1^2+r_2^2}+V$ (see \eqref{relative_eq}) for the Lagrangian 
relative equilibria as a function of $r$.
}
\end{figure}

Figure~\ref{fig:three_body_angularmomentum_energy}a shows the magnitude of the angular momentum $r$ as a function of the  Lagrangian equilateral triangle specified by $b$. One sees that for a given value of $r$, there are two  (or no) Lagrangian equilateral triangle of different size.
The corresponding energies given by the effective potential $V_\text{eff}=\frac12 \frac{r^2}{r_1^2+r_2^2}+V$ at these equilibria are shown in Fig.~\ref{fig:three_body_angularmomentum_energy}b. For a given value of $r$, the smaller Lagrangian triangle has the smaller energy. 

As mentioned above at a relative equilibrium point the body angular momentum vector $\mathbf{J}$ is an eigenvector of $\mathbb{I}^{-1}$.
So, when looking for relative equilibria in general one would like to diagonalize $\mathbb{I}^{-1}$ which is possible when passing to a principal axes frame. We note that the  corresponding
shape coordinates are called \emph{Draght's coordinates} \cite{LittlejohnReinsch97, Iwai87c}. As we were only interested in Lagrangian equilateral triangle configurations in this paper, 
the commonly used Jacobi coordinates were also useful in the study of these relative equilibria since $\mathbb{I}^{-1}$ is diagonal in the third component which corresponds to the direction of the body fixed angular momentum in the 
$xxy$-gauge.

%%%%%%%%%%%%%%%%%%%%%%%%%%%%%%%%%%%%%%
\subsection{Normal form around Lagrangian equilateral triangle relative equilibria}
\label{sec:Normal}

We now apply the procedure explained in Appendix~\ref{algorithm} to compute the Poincar{\'e}-Birkhoff normal form  around the Lagrangian equilibria. We choose unit masses, the parameter $d_0$ in the Morse potential in \eqref{eq:def_Morse_potential} equal to 6 and
the parameter $b$ specifying the side length of the Lagrangian equilateral triangle in \eqref{eq:Jacobi_Lagrange_triangle} equal to 6.5. This gives the Jacobi coordinates
\begin{equation}
\begin{split}
q^e_1 &\equiv r^e_1 =\sqrt{\mu_1}b=\frac{6.5}{\sqrt{2}}\,,  \\
q_2^e &\equiv r^e_2 = \sqrt{\mu_2}b\sqrt{\frac{3}{4}}=\frac{6.5}{\sqrt{2}}\,, \\
q_3^e &\equiv \phi^e =\frac{\pi}{2}\,.
\end{split}
\end{equation}
Solving \eqref{relative_eq} for $r$ we find  $r^e= 19.8302179854%697421410850853495
$. 

The momenta conjugate to the Jacobi coordinates are 
\begin{equation}
\begin{split}
p_1^e &=0\,,\\
p_2^e &=0\,,\\
p_3^3 &=\frac{(r_2^e)^2}{(r_1^e)^2 + (r_2^e)^2}\,
r = 9.9151089927 % 3487107054254267481
\,.
\end{split}
\end{equation}

The eigenvalues of  the matrix $\mathfrak{J}D^2 h_r$ that gives the linearized vector field are 
\begin{equation}
\begin{split}                                           
\pm \ui \, \omega_1 &=  \pm \ui \,  0.2362174000 %21095069619033908429
\,, \\
\pm \ui \, \omega_2 &= \pm \ui \,  0.4693542718 % 45437683812664741998
\,,\\ 
\pm \ui \, \omega_3 &= \pm \ui \,  1.1749259437 % 1197043724436339243
\,, \\                                                        
\pm \ui \, \omega_4 &= \pm \ui \,  1.1984363284 %2141028332592676542
\,.
\end{split}
\end{equation}
So we can immediately read off that the equilibrium is of elliptic linear stability. We note that 
a well established method for determining the  stability of reduced systems is the \emph{reduced energy-momentum method} which was introduced in \cite{Simoetal91}.
For an application to the three-body problem, see also \cite{SchmahStoica06}. The reduced energy-momentum method does however not provide a means to compute higher order normal forms as we will do now following Appdenix~\ref{algorithm}. 

Since we are only interested in demonstrating the basic principle of a normal form computation we will restrict ourselves to the normal of order 4. We start from the fourth order Taylor expansion of the Hamiltonian $h_r$ at the relative equilibrium. It has 212 nonvanishing terms and we refrain from writing them down. The symplectic matrix $M$ which yields the linear symplectic transformation \eqref{eq:M_norm_quadratic_part} after which
the quadratic part of the Hamiltonian assumes the form
\begin{equation}
h^{(2)}_{r,2} = \sum_{k=1}^4 \omega_k (p_k^2 + q_k^2)
\end{equation}
(see \eqref{eq:Hquadratic}) can be defined as
\begin{equation}
M = [ c_1 \text{Re}\, \mathbf{v}_1, c_2 \text{Re}\, \mathbf{v}_2, c_3  \text{Re}\, \mathbf{v}_3, c_4 \text{Re}\, \mathbf{v}_4,    c_1 \text{Im} \, \mathbf{v}_1,  c_2 \text{Im} \, \mathbf{v}_2,  c_3 \text{Im} \, \mathbf{v}_3, c_4 \text{Im} \, \mathbf{v}_4] \label{eq:M_for_three_body}
\end{equation}
where the column vectors are the real and imaginary parts of eigenvectors $\mathbf{v}_k$ of $\mathfrak{J}D^2 h_r$ for the eigenvalues $\ui \, \omega_k$ with coefficients 
\begin{equation}
c_k = \frac{1}{\sqrt{ \text{Re}\, \mathbf{v}_k \cdot  \mathfrak{J}  \text{Im} \, \mathbf{v}_k  }}\,,\quad k =1,\ldots, 4\,.
\end{equation}

\rem{
For these configurations
$\phi=\pi/2$, and we take unit masses, and $b=6.5$, so we find $r_1=r_2=4.59619$. In this case $r=19.83021$ by Eq.~\eqref{relative_eq},
then by \eqref{three_moment} one finds
$p_1=p_2=0$, and $p_3=9.91510$. Then the frequencies are obtained to be
$\omega_1=0.46935 , \omega_2=1.17492 , \omega_3=1.19843$ and $\omega_4=0.23621$.   

$r= 19.8302179854697421410850853495$
position of relative equilibrium $x_0=\sqrt{\mu_1}b=3.25\sqrt{2}=$, $y_0= \sqrt{\mu_2}b\sqrt{\frac{3}{4}}= 1.08333333333333333333333333334\sqrt{18}$, $z_0=\frac{\pi}{2}$
$p_{x0} = 0$, $p_{y0} = 0$, $p_{z0} = 9.91510899273487107054254267481$
} % end rem

We find $M$ 
\begin{equation}
{\tiny
\begin{pmatrix}          
0 & 0.2245619939 & 0 & 0 & 0 & 0 & 0 & 0 \\
0.5952500442 & 0 & 0.5952500442 & -0.6459181965 & 0 & 0 & 0 & 0  \\
-0.5952500442  & 0 & -0.5952500442 & -0.6459181965 & 0  & 0 & 0 & 0 \\
0 & 0 & 0 & 0 & 0.2590186725 & 0 & -0.2590186725 &  0  \\
0 & 0 & 0 & 0 & 0 & 4.4531132913  & 0 & 0 \\
0  & 0 & 0 & 0 & 0.1406084178 & 0 & 0.6993747200 & -0.7740918318  \\
0 & 0 & 0 & 0 & -0.1406084178  & 0 & -0.6993747200  & -0.7740918318 \\
-3.2144619462  & 0 & 0.6462635772  & 0 & 0 & 0 & 0 & 0 
\end{pmatrix}
}\,.
\end{equation}
Following the next step in the procedure described in Appendix~\ref{algorithm} we find that the
normal form of order 4 is given by
\begin{equation}
{\tiny
\begin{split}
h^{(4)}_{r\,\text{NF}} =&  
2.1181531267 % 917529
+ 0.2362174000 % 210951 
\,I_1 + 0.4693542718 %454377
\,I_2  +1.1749259437%1197 
\,I_3 + 1.1984363284% 2141 
\,I_4\\
%%
%(-0,0.2362174000210951) hbar^0 p_1^1 p_2^0 p_3^0 p_4^0 q_1^1 q_2^0 q_3^0 q_4^0 
%(0,0.4693542718454377) hbar^0 p_1^0 p_2^1 p_3^0 p_4^0 q_1^0 q_2^1 q_3^0 q_4^0 
%(-0,1.17492594371197) hbar^0 p_1^0 p_2^0 p_3^1 p_4^0 q_1^0 q_2^0 q_3^1 q_4^0 
%(0,1.19843632842141) hbar^0 p_1^0 p_2^0 p_3^0 p_4^1 q_1^0 q_2^0 q_3^0 q_4^1
%% 
&- 1.4978871558%71867 
\,I_1^2  - 7.7221894156 % 0272 
\,I_1 I_4- 0.9580186364% 034886 
\,I_4^2 + 6.4183166825 % 4332 
\,I_1 I_3 - 8.1361397396 % 40436
\,I_3 I_4 \\ &- 0.8641444715%322564 
\,I_3^2 - 0.2175152611 % 188608 
\,I_1 I_2 - 0.2069751620 % 188574 
\,I_2 I_4 - 0.1815241432 % 409895 
\,I_2 I_3 + 0.0089977794%69467816 
\,I_2^2 \,.
%%
%(1.497887155871867,0) hbar^0 p_1^2 p_2^0 p_3^0 p_4^0 q_1^2 q_2^0 q_3^0 q_4^0 
%(7.72218941560272,0) hbar^0 p_1^1 p_2^0 p_3^0 p_4^1 q_1^1 q_2^0 q_3^0 q_4^1 
%(0.9580186364034886,0) hbar^0 p_1^0 p_2^0 p_3^0 p_4^2 q_1^0 q_2^0 q_3^0 q_4^2 
%(-6.41831668254332,0) hbar^0 p_1^1 p_2^0 p_3^1 p_4^0 q_1^1 q_2^0 q_3^1 q_4^0 
%(8.136139739640436,0) hbar^0 p_1^0 p_2^0 p_3^1 p_4^1 q_1^0 q_2^0 q_3^1 q_4^1 
%(0.8641444715322564,0) hbar^0 p_1^0 p_2^0 p_3^2 p_4^0 q_1^0 q_2^0 q_3^2 q_4^0 
%(0.2175152611188608,0) hbar^0 p_1^1 p_2^1 p_3^0 p_4^0 q_1^1 q_2^1 q_3^0 q_4^0 
%(0.2069751620188574,0) hbar^0 p_1^0 p_2^1 p_3^0 p_4^1 q_1^0 q_2^1 q_3^0 q_4^1 
%(0.1815241432409895,0) hbar^0 p_1^0 p_2^1 p_3^1 p_4^0 q_1^0 q_2^1 q_3^1 q_4^0 
%(0,0.06268608527566027) hbar^0 p_1^1 p_2^2 p_3^0 p_4^0 q_1^0 q_2^0 q_3^1 q_4^0 
%(0,-0.06268608527566027) hbar^0 p_1^0 p_2^0 p_3^1 p_4^0 q_1^1 q_2^2 q_3^0 q_4^0 
%(-0.008997779469467816,0) hbar^0 p_1^0 p_2^2 p_3^0 p_4^0 q_1^0 q_2^2 q_3^0 q_4^0
\end{split} 
} % end tiny
\end{equation}
%\HW{Is there more to say about this result?}
The Hamiltonian $h^{(4)}_{r\,\text{NF}}$ is obtained from the general approach described in this paper.
It yields an integrable nonlinear approximation of the 3-body problem reduced by the non-Abelian symmetry group $SO(3)$ which can be used 
to study the motion in the neighborhood of the Lagrangian equilateral relative equilibria.
Higher order terms can be obtained following the procedure in Appendix~\ref{algorithm}.

%%%%%%%%%%%%%%%%%%%%%%
\subsection{Reconstruction of dynamics}
\label{sec:Reconstruction}

Generally speaking, when a curve $c_\mu$ in the reduced space is given finding the actual curves in the full space is the problem of reconstruction.
This topic is well developed and we refer to \cite{MarsdenMontgomeryRatiu90} for the details. 
We briefly sketch  how the full dynamics can be computed in the three-body reduction case. 

Consider a curve $c_\mu(t)=(r_1(t),r_2(t),\phi(t),p_1(t),p_2(t),p_3(t),u(t),v(t))$ in the reduced space.   Then what is the corresponding curve in the full space of which the projection is $c_\mu$? 
Firstly,
we can find $\mathbf{J}=(J_1,J_2,J_3)$ by \eqref{sphere1}. Then by Eq.~\eqref{angular} it is easy to obtain the angular velocity vector $\xi$ or the corresponding matrix $\Xi$.
After that one has to solve the reconstruction equation \eqref{Omega}
which is not a trivial task because of the non-Abelian structure of $SO(3)$ (see the explanations given in
\cite{MarsdenMontgomeryRatiu90}).
A solution  $g$ gives the $SO(3)$ coordinates, and finally by Eq.~\eqref{rec} their conjugate momenta are obtained.
%Putting all of this  together the actual curve is obtained locally in terms of the identification \eqref{trivialization}.
For a more detailed discussion of the reconstruction in the case of the three-body problem, we refer to  \cite{Montgomery96}.

%%%%%%%%%%%%%%%%%%%%%%
\section{The double spherical pendulum}
\label{sec:Pen}

In this section we study the reduction of the  double spherical pendulum   (for a more 
detailed survey, we refer to \cite{MarsdenScheurle93b}). We again introduce canonical coordinates on the reduced space 
in the light of Sec.~\ref{sec:2}. We use these to compute the Poincar{\'e}-Birkhoff normal form at the relative equilibrium given by a so called  stretched-out solution.

%%%%%%%%%%%%%%%%%%%%%%%%%
\subsection{Reduced equations of motion}

\begin{figure}
\begin{center}
\includegraphics[angle=0,width=7cm]{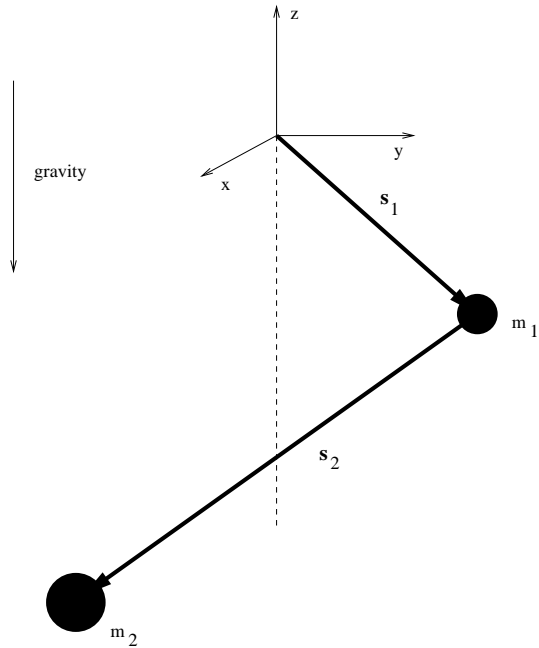}
\end{center}
\caption{\label{fig:double}
The double spherical pendulum.
}
\end{figure}

Consider two coupled spherical pendula with masses $m_1$ and $m_2$ and position vectors
$\mathbf{s}_1$ and $\mathbf{s}_2$ defined as in Fig.~\ref{fig:double},  moving without friction under the influence of a gravitational force $-a\mathbf{k}$ where $a$ is a positive constant and 
$\mathbf{k}$ is the unit vector in the $z$-direction. 

If the lengths of $\mathbf{s}_1$ and $\mathbf{s}_2$ are $l_1$ and $l_2$, respectively,
then the configuration space is $M=S^{2}(l_1) \times S^{2}(l_2)$. The Lagrangian is
\begin{eqnarray}
 L(\mathbf{s}_1,\mathbf{s}_2,\dot{\mathbf{s}}_1,\dot{\mathbf{s}}_2)&=&\frac{1}{2}m_1\,\|\dot{\mathbf{s}}_1\|^2
+\frac{1}{2}m_2\,\|\dot{\mathbf{s}}_1+\dot{\mathbf{s}}_2\|^2 \\
&-&m_1\, a \, \mathbf{s}_1^T\, \mathbf{k}-m_2\, a\,(\mathbf{s}_1+\mathbf{s}_2)^T\, \mathbf{k}\,.
\end{eqnarray}
The system is invariant under  rotations
around the $z$-axis.  So the symmetry group is the Abelian group $S^1$ whose action on $M$ is given by
\begin{equation}
 (\mathbf{s}_1,\mathbf{s}_2)\rightarrow (g_{\theta} \,\mathbf{s}_1,g_{\theta}\,\mathbf{s}_2)\,,
\end{equation}
where $g_{\theta}$ is the rotation by the angle $\theta$ about the $z$-axis. We can identify the Lie algebra of $S^1$ with span$(\mathbf{k})$. An element of the Lie algebra is then an angular velocity vector
of the form  $\omega \mathbf{k}$
with $\omega \in \mathbb{R}$  and the corresponding fundamental vector field is $\omega \, (\mathbf{k} \times \mathbf{s}_1,\mathbf{k} \times \mathbf{s}_2)$. For  the angular momentum, we find according to \eqref{tan_mom_map}
\begin{equation}
 \langle \mathbf{L}(\mathbf{s}_1,\mathbf{s}_2,\dot{\mathbf{s}}_1,\dot{\mathbf{s}}_2), \omega \mathbf{k} \rangle=
\omega (m_1 \dot{\mathbf{s}}_1^T \, (\mathbf{k} \times \mathbf{s}_1) + m_2 (\dot{\mathbf{s}}_1 + \dot{\mathbf{s}}_2)^T \, (\mathbf{k} \times \mathbf{s}_1+ \mathbf{k} \times \mathbf{s}_2))
\end{equation}
or 
 \begin{equation}
 \mathbf{L}= \mathbf{k}^T\, (m_1 (\mathbf{s}_1 \times \dot{\mathbf{s}}_1) + m_2  (\mathbf{s}_1 + \mathbf{s}_2) \times (\dot{\mathbf{s}}_1 +\dot{\mathbf{s}}_2))\, \mathbf{k}\,.
\end{equation} 

If the body frame is chosen such that the $x$-axis coincides with  $\mathbf{s}_1^{\perp}$ and we introduce polar coordinates $(r,\theta)$ in the $xy$-plane then we obtain for
the body-fixed position vectors
\begin{eqnarray}
 \mathbf{r}_1&=&(r_1,0, -\sqrt{l_1^2-r_1^2})\,, \\
 \mathbf{r}_2&=&(r_2 \cos \varphi,r_2 \sin \varphi, -\sqrt{l_2^2-r_2^2})\,.
\end{eqnarray}
Note that through the choice of the sign of the square roots in the last components these equations are restricted to downward pointing configurations.
%\HW{We should first say something about the shape space and its dimension. What actually is the topology? $I\times S^2$ where $I$ is a closed interval? Should we also mention $\mathbb{R}^3_+$ for the three-body case?}
As $r_1, r_2 , \varphi$ are invariant under the group action we can take them as shape coordinates on the three-dimensional shape space $Q=S^{2}(l_1) \times S^{2}(l_2)/S^1$.

The moment of inertia tensor 
\begin{equation}
\mathbb{I}= m_1\, \| \mathbf{r}_1 ^{\perp} \|^2 + m_2\, \|(\mathbf{r}_1+\mathbf{r}_2)^{\perp} \|^2\,,
\end{equation}
where $\mathbf{r}_1^{\perp}$ is the projection of $\mathbf{r}_1$ onto the $xy$-plane,
can be written in terms of the shape coordinates as
\begin{equation}
 \mathbb{I}=(m_1 + m_2) r_1^2 + 2 m_2 r_1 r_2 \cos \varphi +  m_2 r_2^2\,.
\end{equation}
Accordingly, we get for the mechanical connection
\begin{equation}
 A_{r_1}=-\frac{m_1 m_2 r_2 \sin \varphi}{\mathbb{I}},\hspace{2mm} A_{r_2}=\frac{m_1 m_2 r_1 \sin \varphi}{\mathbb{I}},
\hspace{2mm} A_{\varphi}=\frac{m_1 m_2 r_2 (r_1 \cos \varphi +r_2)}{\mathbb{I}},
\end{equation} 
and the entries of the matrix $d$ which gives the horizontal metric are
\begin{eqnarray*}
 d_{11}&=&\frac{l_1^2 (m_1+m_2)}{2 (l_1^2-r_1^2)}
-\frac{m_1^2 m_2^2 r_2^2 \sin^2 \varphi}{(m_1+m_2) r_1^2+
2 m_2 r_1 r_2 \cos \varphi + m_2 r_2^2}, \\
d_{12}&=& \frac{1}{2} m_2 \left(\cos \varphi + r_1 r_2 \left(\frac{1}{\sqrt{l_1^2
-r_1^2} \sqrt{l_2^2- r_2^2}}+\frac{2 m_1^2 m_2 sin^2 \varphi}{(m_1+m_2)
 r_1^2 +2 m_2 r_1 r_2 \cos \varphi + m_2 r_2^2}\right)\right), \\ 
d_{13} &=& \frac{1}{2} m_2 r_2 \left(-1+\frac{2 m_1^2 m_2 r_2
(r_1 \cos \varphi+ r_2)}{(m_1+m_2) r_1^2
+2 m_2 r_1 r_2 \cos \varphi+ m_2 r_2^2}\right) \sin \varphi, \\
d_{22}&=& m_2 \left(\frac{l_2^2}{2 l_2^2-2 r_2^2}-\frac{m_1^2 m_2 r_1^2 \sin^2 \varphi}
{(m_1+ m_2) r_1^2+2 m_2 r_1 r_2 \cos \varphi+m_2 r_2^2}\right), \\
d_{23}&=& -\frac{m_1^2 m_2^2 r_1 r_2 (r_1 \cos \varphi+ r_2 \sin \varphi)}
{(m_1+m_2) r_1^2+2 m_2 r_1 r_2 \cos \varphi+ m_2 r_2^2}, \\
d_{33}&=& \frac{1}{2} m_2 r_2^2 \left(1-\frac{2 m_1^2 m_2 (r_1 \cos \varphi+ r_2^2}{(m_1+m_2) r_1^2+2 m_2
r_1 r_2 \cos \varphi+ m_2 r_2^2}\right)\,.
\end{eqnarray*}

%%%%%%%%%%%%%%%%%%%%%%

The left-action of the group on the tangent bundle is trivial, so $\mathbf{J}=\mathbf{L}$.
The conjugate momenta of the shape coordinates are given by
\begin{equation}
 p_\alpha =\frac{\partial L}{\partial q_\alpha}=\sum^{3}_{\beta=1}g_{\alpha \beta} \dot{q}_\beta + \mathbf{J} \, \mathbf{A}_\alpha. \label{spher_moment}
\end{equation}
Thus the Hamiltonian takes the form
\begin{equation}
h=\frac{1}{2} \mathbb{I}^{-1} \mathbf{J}^2+
\frac{1}{2} \sum^{3}_{\alpha,\beta=1} d^{\alpha \beta} (p_{\alpha}-\mathbf{J} \, \mathbf{A}_\alpha)(p_\beta-\mathbf{J} \, \mathbf{A}_\beta)+V, 
\end{equation}
where $V=-m_1 \,a\,\sqrt{l_1^2-r_1^2}-m_2 \,a\,(\sqrt{l_1^2-r_1^2}+\sqrt{l_2^2-r_2^2})$ is the potential. Observe here that as $\mathbf{L}$ is conserved,
$\mathbf{J}$ is also conserved and can be viewed as a parameter.

The reduced equations of motion are 
\begin{equation}
\begin{split}
 \dot{q}_{\alpha}=&\frac{\partial h_r}{\partial p_{\alpha}}=d^{\alpha\beta}
 (p_{\beta}-  \mathbf{J} \, A_{\beta} ),\\
 \dot{p_{\alpha}}=&-\frac{\partial h_r}{\partial q_{\alpha}}=-\frac{\partial}{\partial q_\alpha}(\frac{1}{2} \{ \mathbb{I}^{-1} \mathbf{J}^2 +
 d^{\alpha\beta}(p_{\alpha}- \mathbf{J} \, A_{\alpha})\, (p_{\beta}- \mathbf{J} \, A_{\beta} )\}
+ V(q))\,,
\end{split}
\end{equation}
where we denote by $r$ the $z$-component of the conserved angular momentum $\mathbf{J}$.

%%%%%%%%%%%%%%%%%%%%%%%%%%%%%
\subsection{Relative equilibria}
\label{sec:dsp_rel_equilibria}

\begin{figure}
\begin{center}
\includegraphics[angle=0,width=7cm]{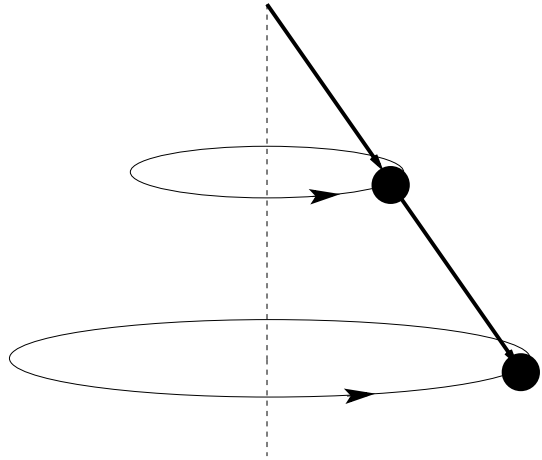}
\end{center}
\caption{\label{fig:streched}
A stretched-out relative equilibrium solution of the double spherical pendulum in which the two masses are aligned with the point of suspension and move along circles. 
}
\end{figure}

\begin{figure}
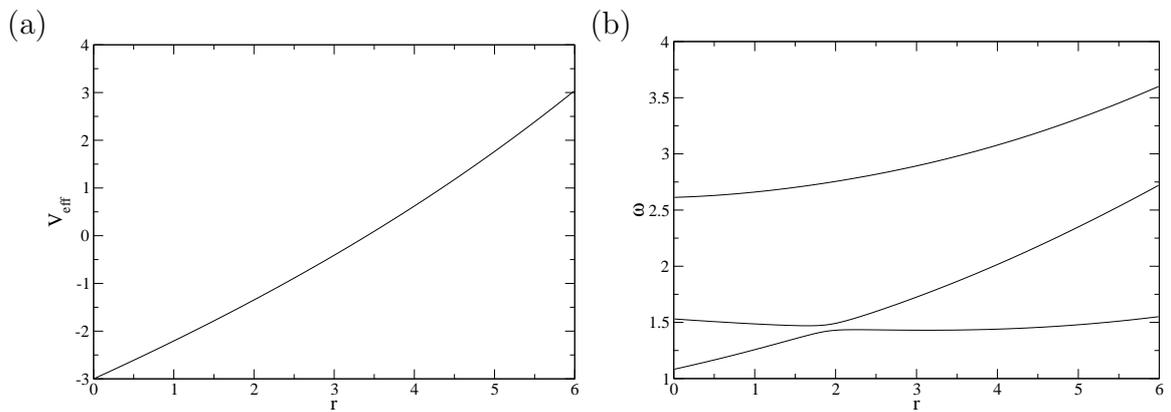

\begin{center}
\raisebox{5cm}{(a)}\includegraphics[angle=0,width=7cm]{double_pendulum_energy} 
\raisebox{5cm}{(b)}\includegraphics[angle=0,width=7cm]{double_pendulum_frequencies}
\end{center}
\caption{\label{fig:dsp}
(a) Energy of the stretched-out relative equilibrium as given by the effective or amended potential $V_{\text{eff}} = \frac{1}{2} \mathbb{I}^{-1} \mathbf{J}^2 + V(q)$ (see \eqref{amended}) as a function of the angular momentum $r$.
(b) Frequencies of the stretched-out relative equilibrium as a function of $r$.}
\end{figure}

Because of the triviality of the coadjoint action 
the conditions to have a  relative equilibrium \eqref{eq:cond_rel_equil_gen} reduce to
\begin{equation}
\begin{split}
p_{\alpha} =  \mathbf{J} \, A_{\alpha} \,, \\
\frac{\partial}{\partial q_\alpha }\left( \frac{1}{2} \mathbb{I}^{-1} \mathbf{J}^2 + V(q) \right) &=0\,. \label{amended}
\end{split}
\end{equation}
As shown in  \cite{MarsdenScheurle93b} there are two types of relative equilibria: the so called cowboy branch and the stretched-out solution that we will concentrate on in the following and which is shown in Fig.~\ref{fig:streched}.
For a stretched-out relative equilibrium, we have $\varphi=0$. We find $r_1$ and $r_2$ from solving \eqref{amended} (using the computer algebra program Maple). The corresponding momenta are obtained from \eqref{spher_moment}. We will in the following choose all  parameters to have unit values, i.e.
$m_1=l_1=m_2=l_2=a=1$. The energy of the stretched-out relative equilibrium as a function of $\mathbf{J}$ for this choice of parameters is shown in Fig.~\ref{fig:dsp}a.

%%%%%%%%%%%%%%%%%%%%%%
\subsection{Normal form around stretched-out relative equilibria}
\label{stretched}

%The relative equilibria for the double spherical pendulum were studied
%in \cite{MarsdenScheurle93b} so we do not find them and do a stability test. 
%As an example we consider the stretched out solution (see Fig.~\ref{fig:streched}).

The stretched-out relative equilibria are known to be stable \cite{MarsdenScheurle93b}. In agreement with this result we find that the eigenvalues of the matrix $\mathfrak{J}D^2 h_r$ associated with the linearized vector field of the reduced system at the stretched-out relative equilibria
are purely imaginary. The frequencies $\omega_k$, $k=1,2,3$, are shown in Fig.~\ref{fig:dsp}b as a function of $r$.
 
For the normal form computation, we consider the relative equilibrium point which has $r=1$. 
For the position of this relative equilibrium, we find
\begin{equation}
\begin{split}
q_1^e&\equiv r^e_1 = 0.4425598655\,, %30592352522120140694 
\\  q_2^e &\equiv r^e_2 = 0.5656579210%71478730859028380941 
\\  q_3^e &\equiv  \varphi^e = 0\,,\\ 
                             p_1 &=  0\,, \\
                              p_2 &=  0\,, \\
         p_3 &=       0.4704091824 %98242938678237166424
         \,.
\end{split}
\end{equation}
The Taylor expansion of the reduced Hamiltonian  to order 4 has 186 nonvanishing terms at this relative equilbrium.
The eigenvalues associated with the linearized vector field are 
\begin{equation}
\begin{split}      
\pm \ui \,\omega_1  &= \pm \ui\,  1.2572610531 %23528
\\                                             
 \pm \ui\, \omega_2  &= \pm \ui\, 1.4864684140%11076
 \\                                                           
 \pm \ui \,\omega_3  &= \pm \ui \, 2.6603546311% 09321 
\end{split}
\end{equation}
We define the symplectic matrix $M$ which transforms the quadratic part of the Hamiltonian to the form $\sum_{k=1}^3 \omega_k(p_k^2+q_k^2)$ analogously to \eqref{eq:M_for_three_body}  
 in  Sec.~\ref{sec:Normal}. 
 We find
\begin{equation} \begin{tiny}
M=
\begin{pmatrix}
0.3720476175 & -0.3116300067 & 0.4202281883 & 0 & 0 & 0\\
-0.5712604029 & -0.3238407419 & -0.5369033436 & 0 & 0 & 0 \\
0 & 0 & 0 & 2.2390869882 & -0.2974563002 & -2.2029538091 \\
0 & 0 & 0 & -0.8426837965 & -1.6576611267 &  1.8964495311 \\
0 & 0 & 0 &  0.3873035344 & -1.4365053634 & -1.4081719163 \\
-0.6854444155 & 0.0612664420 & -0.2510237658 & 0 & 0 & 0
\end{pmatrix} \,.
\end{tiny}
\end{equation}
Following the procedure described in Appendix~\ref{algorithm}
we get for the $4$th order normal form
\begin{equation}
\begin{split}
h_{r\,\text{NF}}^{(4)} =&
-2.2056999577% 853963 
+ 1.2572610531%23528 
I_1 + 1.4864684140%11076 
I_2 + 2.6603546311%09321 
I_3 \\
%(0,1.257261053123528) hbar^0 p_1^1 p_2^0 p_3^0 q_1^1 q_2^0 q_3^0 
%(0,1.486468414011076) hbar^0 p_1^0 p_2^1 p_3^0 q_1^0 q_2^1 q_3^0 
%(0,2.660354631109321) hbar^0 p_1^0 p_2^0 p_3^1 q_1^0 q_2^0 q_3^1 
%%
&+ 0.0467015469%467138 
I_1^2 - 5.8213832524%3255 
I_1 I_3 + 0.0875786340%0963339 
I_3^2 +  0.1772800788%017088 
 I_1 I_2\\& - 0.0932948515%85636 
I_2 I_3 - 0.0419637147%8522198 
I_2^2\,.
%%
%(-0.0467015469467138,0) hbar^0 p_1^2 p_2^0 p_3^0 q_1^2 q_2^0 q_3^0 
%(5.82138325243255,0) hbar^0 p_1^1 p_2^0 p_3^1 q_1^1 q_2^0 q_3^1 
%(-0.08757863400963339,0) hbar^0 p_1^0 p_2^0 p_3^2 q_1^0 q_2^0 q_3^2 
%(-0.1772800788017088,0) hbar^0 p_1^1 p_2^1 p_3^0 q_1^1 q_2^1 q_3^0 
%(0.093294851585636,0) hbar^0 p_1^0 p_2^1 p_3^1 q_1^0 q_2^1 q_3^1 
%(0.04196371478522198,0) hbar^0 p_1^0 p_2^2 p_3^0 q_1^0 q_2^2 q_3^0 
\end{split}
\end{equation}
%HW{More comments on this result?}
The Hamiltonian $h^{(4)}_{r\,\text{NF}}$ yields an integrable nonlinear approximation of the 
double spherical pendulum reduced by the Abelian symmetry group $SO(2)$ in the neighborhood of the stretched-out relative equilibrium.
Similarly to the 3-body case in Sec.~\ref{sec:Normal} the normal form of the reduced system is obtained from the general approach in this paper which demonstrates the effectiveness and generality of the approach.

%%%%%%%%%%%%%%%%%%%%%%
\section{Conclusions and Outlook }
\label{sec:conclusions}

%
% what we did
%
In this paper we provided a general perspective on the construction of canonical coordinates for the reduced spaces of Hamiltonian systems given by cotangent bundles with a free Lie group action. The general approach presented in this paper allows one to treat the reduction of Abelian and non-Abelian group actions on the same footing.  The case of simple mechanical systems was studied in detail. The approach was illustrated for a 3-body problem and the double spherical pendulum which involve non-Abelian and Abelian symmetries, respectively. We used the  canonical coordinates to compute the  Poincar{\'e}-Birkhoff normal forms at the relative equilibria given by the Lagrangian equilateral triangle configuration  in the 3-body problem and the stretched-out solution of the double spherical pendulum.  The Poincar{\'e}-Birkhoff normal form gives a nonlinear approximation of the local dynamics of the reduced system in the neighborhood of the relative equilibria. This goes beyond the well established reduced energy-momentum method which only give the linear stability of the relative equilibria  \cite{Simoetal91}, and enables one, e.g., to give nonlinear approximations of the center manifolds of relative equilibria. The use of  a Poincar{\'e}-Birkhoff normal form for the computation of the center manifolds of  saddle type equilibria has in recent years been demonstrated in the study of reaction type dynamics \cite{Uzeretal02,WaalkensBurbanksWigginsb04}. The study of this paper allows one to carry over these results to the case of saddle type relative equilbria which induce reaction type dynamics in rotating molecules \cite{CiftciWaalkens11}. 
In this context also the reconstruction of the full dynamics from the reduced one which we illustrated for the example of the 3-body problem is of importance.
 
%
% isotropy 
%
In this paper we excluded Lie group actions with isotropy  which have been studied, e.g.,  in \cite{RobertsScmahStoica06} or \cite{KozinRobertsTennyson00} and \cite{IwaiYamaoka05}
for the case of  $n$-body systems. Our future studies concern how isotropy can be incorporated  in the approach presented in this paper to do, e.g,  a normal form analysis for
3-body systems with linear equilibrium configurations.
%
% gauge independence
%
Another related problem is the development of a normal form algorithm which is coordinate independent or in the jargon of \cite{LittlejohnReinsch97}  \emph{gauge independent}.
This problem is considered in \cite{LittlejohnMitchell03} and its generalizations to cotangent bundles in general seems worth studying.

%
% nonholonimic
%
Recently considerable progress was made in nonholonomic mechanics (see, e.g.,  \cite{Blochetal96, Blochetal09} and the references therein). It might be possible to use these techniques to develop a non-canonical Poincar{\'e}-Birkhoff normal form related  the results of this paper.

%
%
% quantum normal form
%
Finally we mention that it would be interesting to transfer the results of this paper  to
 quantum mechanical systems.
Analogously to the Poincar{\'e}-Birkhoff normal form of an equilibrium point there is a quantum normal form built on the symbol calculus of pseudo differential operators by which one can locally approximate a quantum Hamilton operator.  
In the case of elliptic equilibria this allows one to compute quantum energy spectra with high precision. For saddle type equilibria, the quantum normal form can be used 
to compute efficiently quantum reaction rates and the associated Gamov-Siegert resonances (see \cite{SchubertWaalkensWiggins06,Waalkensetal08} and also the references therein for quantum normal forms in general).
It would be interesting to transfer these results to relative equilibria of rotational symmetry reduced molecular systems.
The dependence of the quantization of quantum reaction rates in the hydrogen exchange 
reaction as a function of the angular momentum has, e.g., been studied in \cite{ChatfieldMielkeAllisonTruhlar2000} using  \emph{ab initio}
quantum computations. It would be interesting to compare these results to a quantum normal computation. 
The geometric approach for quantum 3-body problems  presented in  \cite{Iwai87c} could be very useful for this purpose.

%%%%%%%%%%%%%%%%%%%%%%
\appendix
%\label{sec:appendix}

%%%%%%%%%%%%%%%%%%%%%%%%%%%%%%%
\section{Poisson brackets in an anholonomic frame}
\label{anholonomic}
We give the Poisson brackets in an anholonomic frame,
the details can be found in \cite{LittlejohnReinsch97}.
Let $M$ be a manifold and $s_1,...,s_n$ be a local coordinate system on $M$.
Consider a regular Lagrangian $L:TM \rightarrow \mathbb{R}$ which then is a function of
the coordinates $s_1,...,s_n,\dot{s}_1,...,\dot{s}_n$, and a local frame $X_1,...,X_n$ 
on $M$ with
\begin{equation}
 X_i= \sum_{j=1}^{n} a_i^j \partial_j\,,
\end{equation}
where $\partial_i = \partial / \partial s_i$.  
From
\begin{equation}
\dot{s}= \sum_{j=1}^{n} \dot{s}_i \partial_i = \sum_{j=1}^{n} v_i X_i,
\end{equation}
we obtain  that
\begin{equation}
 \frac{\partial \dot{s}_i}{\partial v_j} = a_j^i\,.   \label{partial}
\end{equation}
Let $\bar{L}(s_i,v_i)=L(s_i,\dot{s}_i)$. Then the conjugate momenta
$p_i=\frac{\partial L}{\partial \dot{s}_i}$ and 
$\pi_i=\frac{\partial \bar{L}}{\partial v_i}$ are related by
\begin{equation}
\pi _i = a_i^j p_j \,. \label{pi} 
\end{equation}
If the Poisson bracket of two functions $f$ and $g$ on $M$ in terms of the coordinates $s_1,...,s_n,p_1,...,p_n$  has the canonical form
\begin{equation}\label{pbe}
 \{f,g\}=\sum_{i=1}^{n} \frac{\partial f}{\partial s_i} \frac{\partial g}{\partial p_i}
-\frac{\partial f}{\partial p_i} \frac{\partial g}{\partial s_i}\,,
\end{equation}
then using  Equations~ \eqref{partial} and \eqref{pi} the Poisson brackets of the coordinates $s_1,...,s_n,\pi_1,...,\pi_n$ are given by 
\begin{equation}
 \{ s_i , s_j \} =0, \quad  \{ s_i,\pi_j\}=a_i^j, \quad 
\{\pi_i,\pi_j\}=-\sum_{k=1}^{n} c_{ij}^k \pi_k, \label{eq:PoissonAnholonFrame}
\end{equation}
where $[ X_i,X_j ] =  \sum _{k=1}^{n} c_{ij}^k X_k$. Then the
Poisson bracket (\ref{pbe}) becomes
\begin{equation}
 \{f,g\}=\sum_{i=1}^{n} a_j^i \, \left( \frac{\partial f}{\partial s_i} \frac{\partial g}{\partial \pi_j}
-\frac{\partial f}{\partial \pi_j} \frac{\partial g}{\partial s_i} \right) -c_{ij}^k \, \pi_k
\frac{\partial f}{\partial \pi_j}\, \frac{\partial g}{\partial \pi_i}  \,.
\end{equation}

\section{Lie-Poisson structures}
\label{lie-poisson}
We recall some basics of Lie Poisson structures on Lie groups.
For the details we refer to \cite{Marsden92}.
Let $G$ be a Lie group and $\mathfrak{g}$ be its
Lie algebra with the Lie bracket $[\,,\,]:\mathfrak{g} \times \mathfrak{g}
\rightarrow \mathfrak{g}$. The dual space $\mathfrak{g}^*$ of $\mathfrak{g} $
is a Poisson manifold with either of the following brackets
\begin{equation}
\{f,k\}_{\pm}(\mu)=\pm \left< \mu \, ,\, \left[ \frac{\delta f}{\delta \mu}\,,\,\frac{\delta k}
{\delta \mu} \right] \right>.
\end{equation}
Here $\mathfrak{g}^{**}$ is identified with
$\mathfrak{g}$ in the sense that $\delta f / \delta \mu \in \mathfrak{g}$
is defined by
$\left< \nu \, , \, \delta f / \delta \mu \right> =\mathbf{D}f(\mu)$,
where $\mathbf{D}$ is the derivative. Let $B=e_1,...,e_m$ be a basis of
$\mathfrak{g}$ and $e_1^*,...,e_m^*$ be the dual basis of $B$ in
$\mathfrak{g}^*$ with the corresponding coordinates
$\xi_1,...,\xi_m$ on $\mathfrak{g}$ and
$\mu_1,...,\mu_m$ on $\mathfrak{g}^*$, respectively.
Then the Lie Poisson bracket is given by
\begin{equation}
\{f,k\}_{\pm}=\pm \mu_a \, \gamma_{bc}^a \, \frac{\delta f}{\delta \mu}_b \,
\frac{\delta f}{\delta \mu}_c,
\end{equation}
where $\left[ e_a \, , \, e_b \right]=\gamma_{ab}^c \, e_c$.

Let $H: \mathfrak{g}^* \rightarrow \mathbb{R}$ be a Hamiltoninan
function. Then the general Lie Poisson
equations determined by
$\dot{F}=\{F,H\}$ read
\begin{equation}
\dot{\mu}=\pm \mbox{ad}^*_{\delta H / \delta \mu} \mu,
\end{equation}
where $\mbox{ad}_{\xi}:\mathfrak{g} \rightarrow \mathfrak{g}$
is the linear map $\zeta \rightarrow \left[\xi \,,\,\zeta \right]$, and
$\mbox{ad}_{\xi}^*: \mathfrak{g}^* \rightarrow \mathfrak{g}^*$ is its dual.

%%%%%%%%%%%%%%%%%%%%%%%%%%%%%%%
\section{The Poincar{\'e}-Birkhoff normal form}
\label{algorithm}

Let $H_2$ denote the quadratic Hamiltonian which gives the linearized Hamiltonian vector field at the equilibrium.
One says that a Hamiltonian $H$ is in normal form if
$H$ Poisson commutes with its quadratic part, i.e.
\begin{equation}
\{ H_2 , H \} := \sum_{k=1}^\dof   \big( \frac{\partial H_2}{\partial q_k} \frac{\partial H}{\partial p_k}  -  \frac{\partial H}{\partial q_k} \frac{\partial H_2}{\partial p_k} \big)=0 \,,
\end{equation}
where $f$ is the number of degrees of freedom.
In general $H$ is not in normal form. However, for any given order $n_0$ of the Taylor expansion of $H$ one can find a symplectic transformation to new phase space coordinates
 in terms of  which the transformed $H$ truncated at order $n_0$ is in normal form. 
This symplectic transformation is constructed from a sequence of the form \label{Order $(q,p)$ or $(p,q)$ should be made consistent in the whole paper.}
%which changes the phase space coordinates $\vz \equiv (\vq, \vp)$ order by order according to
\begin{equation}
(q_i,p_i)\equiv \vz\equiv \vz^{(0)} \mapsto \vz^{(1)} \mapsto \vz^{(2)} \mapsto \vz^{(3)}
\mapsto \ldots \mapsto \vz^{(n_0)}\,,
\label{eq:seq_trafos_z}
\end{equation}
where  $\vz^{(n)} $ is obtained from $\vz^{(n-1)} $ 
by means of a symplectic transformation 
\begin{equation}
\vz^{(n-1)} \mapsto \vz^{(n)} =  \phi_{W_n} \vz^{(n-1)} 
\label{eq:zn-1twozn}
\end{equation}
generated by a polynomial $W_n(\vz)$ of order $n$, i.e.
\begin{eqnarray}
W_n \in  \mathcal{W}^n := \mathrm{span} \left\{ q_1^{\alpha_1} \ldots q_\dof^{\alpha_\dof}
    p_1^{\beta_1} \ldots p_\dof^{\beta_\dof}  :   |\alpha|+|\beta|= n \right\} \, .
\label{2-09_cl}
\end{eqnarray}
Here  $|\alpha| = \sum_{k=1}^\dof \alpha_k$, $|\beta| = \sum_{k=1}^\dof \beta_k$.
More precisely, the  $\phi_{W_n}$ in \eqref{eq:zn-1twozn} denote the time-one maps of the flows generated by the Hamiltonian vector fields corresponding to the polynomials $W_n$ (see \cite{Waalkensetal08} for the details). 
The maximum order $n_0$ in \eqref{eq:seq_trafos_z} is the desired order of accuracy at which the expansion will be terminated and truncated. 

Expressing the Hamiltonian $H$ in the coordinates  $\vz^{(n)}$, $n=1,\ldots,n_0$, 
we get a sequence of Hamiltonians $H^{(n)}$,
\begin{equation}
H\equiv H^{(0)} \rightarrow H^{(1)} \rightarrow H^{(2)} \rightarrow H^{(3)}
\rightarrow \ldots \rightarrow H^{(n_0)}\,, 
\label{eq:seq_trafos_cl}
\end{equation}
where for $n=1,\ldots,n_0$,
$
H^{(n)} (\vz^{(n)}) = H^{(n-1)} (\vz^{(n-1)}  ) = H^{(n-1)} ( \phi_{W_n}^{-1} \vz^{(n)})     
$, i.e. 
\begin{equation}
H^{(n)}  = H^{(n-1)} \circ  \phi_{W_n}^{-1}\,.
\end{equation}
To avoid a proliferation of notation we will in the following neglect the superscripts $(n)$ for the phase space coordinates.

In the first transformation in \eqref{eq:seq_trafos_z} we
shift the equilibrium point
$\vz_0$ to the origin, i.e. $\vz \mapsto  \phi_{W_1} (\vz) := \vz - \vz_0$. This gives  
\begin{equation}
  H^{(1)}(\vz) = H^{(0)} (\vz+\vz_0) \, .
\label{2-02_cl}
\end{equation}
The next steps of the normal form procedure rely on  the power series
expansions of $H^{(n)}$, %\footnote{probably we want 1 instead of $n$ here since we actually do not Taylor expand the $H^{(n)}$ for $n\ge 3$.}
\begin{equation}
  H^{(n)}(\vz) = E_0 + \sum_{s=2}^\infty H_s^{(n)}(\vz) \, ,
\label{2-03_cl}
\end{equation}
where the $H_s^{(n)}$ are homogenous polynomials in $ \mathcal{W}^n $:
\begin{equation}
  H_s^{(n)}(\vz)  =  \displaystyle \sum_{|\alpha|+|\beta|=s}
    \frac{H_{\alpha_1,\ldots,\alpha_\dof,\beta_1,\ldots,\beta_\dof}^{(n)}}{\alpha_{1}!
      \ldots \alpha_{\dof}! \beta_{1}! \ldots \beta_{\dof}! }   
    \, q_1^{\alpha_1} \ldots q_\dof^{\alpha_\dof} p_1^{\beta_1} \ldots
    p_\dof^{\beta_\dof}  \, .
 \label{2-04_cl}
\end{equation}
For $n=1$, the coefficients in \eqref{2-04_cl} are given by the Taylor expansion of $H^{(1)}$ at the origin
\begin{equation}
   H_{\alpha_1,\ldots,\alpha_\dof,\beta_1,\ldots,\beta_\dof}^{(1)} 
  = \displaystyle \left.
  \prod_{k,l=1}^\dof \frac{\partial^{\alpha_k}}{\partial
    q_k^{\alpha_k}} \frac{\partial^{\beta_l}}{\partial
    p_l^{\beta_l}} 
  H^{(1)}(\vz) \right|_{\vz={\bf 0}} \!\!\! .
  \label{2-04.1_cl}
\end{equation}
For $n\ge 3$, the coefficients in \eqref{2-04_cl}   are obtained recursively. For $n=2$, i.e. 
the second step in the sequence of transformations \eqref{eq:seq_trafos_z}, the coefficients in  \eqref{2-04_cl} are determined by
 a linear transformation of the phase space coordinates according to
\begin{equation}
\vz \mapsto  \phi_{W_2} (\vz) :=  M\, \vz\,. \label{eq:M_norm_quadratic_part}
\end{equation}
Here,  $M$ is a symplectic $2f \times 2f$ matrix which is chosen in such a way 
 that
the quadratic part of the transformed  Hamiltonian function
\begin{equation}
  H^{(2)}(\vz) = H^{(1)}(M^{-1} \vz)
\label{2-04.2}
\end{equation}
assumes the form
\begin{equation}
 H^{(2)}_2(q,p)=\displaystyle \sum_{k=1}^{f} \frac{\omega_k}{2}(p_k^2+q_k^2).  \label{eq:Hquadratic}
\end{equation}

For the first two steps in the sequence \eqref{eq:seq_trafos_z}, we actually did not give explicit expressions
for the generating functions $W_1$ and $W_2$. For conceptual reasons (and to justify the notation) it is worth
mentioning that such expression can be determined (see \cite{Waalkensetal08}).  The next  steps in  \eqref{eq:seq_trafos_z}
though rely on the explicit computation of the generating functions $W_n$ with $n\ge3$.
To this end it is convenient to introduce the adjoint operator associated with a phase space function $A$:
\begin{equation}
  \mathrm{ad}_A : B \mapsto \mathrm{ad}_A B \equiv \{ A, B \}
  \, .
\label{2-08_cl}
\end{equation}
The transformation \eqref{eq:zn-1twozn} then
leads to a transformation of the Hamilton function $H^{(n-1)}$ to $H^{(n)}$ with $n \ge 3$ which in terms of the adjoint operator 
 reads
\begin{equation}
  H^{(n)} = \sum_{k=0}^\infty \frac{1}{k!} \left[ \mathrm{ad}_{W_n} \right]^k
  H^{(n-1)} \, .
\label{2-10_cl}
\end{equation}
In terms of the Taylor expansion defined in
Eqs.~\eqref{2-03_cl}-\eqref{2-04.1_cl} the transformation introduced by
Eq.~\eqref{2-10_cl} reads
\begin{equation}
  H^{(n)}_s = \sum_{k=0}^{\left\lfloor \frac{s}{n-2} \right\rfloor}
  \frac{1}{k!} \left[ \mathrm{ad}_{W_n} \right]^k H^{(n-1)}_{s-k(n-2)} \, ,
\label{2-11_cl}
\end{equation}
where $\lfloor \cdot \rfloor$ gives the integer part of a number,
i.e., the `floor'-function. 

Using Eq.~\eqref{2-11_cl} one finds that
the transformation defined by \eqref{2-10_cl}  satisfies the following
important properties for $n \ge 3$. Firstly, at step $n$, $n\ge3$, the terms of order less than $n$ in
the power series of the Hamiltonian are unchanged, i.e.
\begin{equation}
  H_s^{(n)} = H_s^{(n-1)} \, , \;\;\; \mathrm{for} \;\;\; s < n \, ,
\label{2-12_cl}
\end{equation}
so that, in particular, $H_2^{(n)} = H_2^{(2)}$.
Defining
\begin{equation}
  \mathcal{D} \equiv \mathrm{ad}_{H_2^{(2)}} = \{ H_2^{(2)} , \cdot \} 
\label{2-14_cl}
\end{equation}
we get for the term of order $n$,
\begin{equation}
  H_n^{(n)} = H_n^{(n-1)} - \mathcal{D} W_n \, .
\label{2-13_cl}
\end{equation}
This is the so-called {\it homological  equation} which will determine the generating functions $W_n$ for $n\ge3$ from 
requiring   $\mathcal{D}
H_n^{(n)} = 0$, or equivalently $H_n^{(n)}$ to be in the kernel of the
restriction of $\mathcal{D}$ to $\mathcal{W}^n$. 
In view of \eqref{2-13_cl}  this condition yields
\begin{equation}
  H_n^{(n-1)} - \mathcal{D} W_n \in \mathrm{Ker} \,\mathcal{D} |_{\mathcal{W}^n} \, .
\label{2-15_cl}
\end{equation}
Section 3.4.1 of Ref.~\cite{Waalkensetal08} provides the explicit
procedure of finding the solution of Eq.~\eqref{2-15_cl}.  In the generic situation where the
linear frequencies $\omega_1,\ldots,\omega_\dof$ in \eqref{eq:Hquadratic} are
rationally independent, i.e. $m_1\omega_1+\ldots+m_\dof \omega_\dof =0$
implies $ m_1=\ldots=m_\dof=0$ for all integers $m_1,\ldots,m_\dof$, it
follows that for odd $n$, $H_n^{(n)} =0$, and for even $n$,
\begin{equation}
  H_n^{(n)} \in \mathrm{span} \left\{  I_1^{\alpha_1} I_2^{\alpha_2} I_3^{\alpha_3}
    \ldots I_\dof^{\alpha_\dof}  : |\alpha|=n/2 \right\} \, ,
\label{2-15.2}
\end{equation}
where $I_k = (q_k^2 + p_k^2)/2$, with $k = 1, \ldots, \dof$.

Applying the transformation \eqref{2-10_cl}, with the generating function
defined by \eqref{2-13_cl}, for $n=3,\ldots,n_0$, and truncating the
resulting power series at order $n_0$ one
arrives at the Hamiltonian $H_{\mathrm{NF}}^{(n_0)}$ corresponding to
the $n_0^{\mathrm{th}}$ order {\it normal form} (NF) of the
Hamiltonian $H$:
\begin{equation}
  H_{\mathrm{NF}}^{(n_0)}(\vz) = E_0 + \sum_{s=2}^{n_0} H_s^{(n_0)}(\vz) \, .
\label{2-15.5_cl}
\end{equation}
The normalized Hamiltonian $H_{\mathrm{NF}}^{(n_0)}$ is an
$n_0^\mathrm{th}$ order approximation of the original Hamiltonian $H$ obtained from
expressing $H$ in terms of the \emph{normal form coordinates} $\vz_{\mathrm{NF}} $ which in turn are obtained from the symplectic transformation of the original coordinates $\vz=(q_i,p_i)$ 
\begin{equation}
\vz_{\mathrm{NF}} = \phi (\vz) =  (\phi_{W_{n_0}} \circ \phi_{W_{n_0-1}}    \circ   \cdots  \circ \phi_{W_{2}}    \circ  \phi_{W_{1}})(\vz)  \,.
\label{eq:def_U}
\end{equation} 
This  is why one can use
$H^{(n_0)}_{\mathrm{NF}}$ instead of $H$ to analyze the dynamics in the neighborhood of an equilibrium. 

The procedure to compute $H_{\mathrm{NF}}^{(n_0)} $ and the corresponding coordinate transformation is algebraic
in nature, and can be implemented on a computer. A computer program is freely available from \cite{Software}.

%%%%%%%%%%%%%%%%%%%%%%

%\bibliographystyle{unsrt}
%\bibliography{redbib}

\end{document}